\documentclass[12pt,twoside,reqno]{amsart}

\usepackage{bm}

\usepackage{amssymb,amsmath,amstext,amsthm,amsfonts,amscd,xcolor}

\usepackage{bbm}
\usepackage{dsfont}

\usepackage{amsthm, enumerate}

\usepackage[ansinew]{inputenc}
\usepackage{graphicx}
\usepackage[mathscr]{eucal}
\usepackage{hyperref}
\usepackage{tikz-cd}

\newcommand{\R}{\mathbb{R}}

\newcommand{\N}{\mathbb{N}}
\newcommand{\Z}{\mathbb{Z}}

\newcommand{\GL}{{\rm GL}}

\newcommand{\Mat}{{\rm Mat}}

\newcommand{\Pp}{\mathbb{P}}

\newcommand{\EE}{\mathbb{E}}
\newcommand{\Escr}{\mathscr{E}}

\newcommand{\cocycles}{\mathcal{C}}

\theoremstyle{plain}
\newtheorem{theorem}{Theorem}[section]
\newtheorem{proposition}{Proposition}[section]
\newtheorem{corollary}[proposition]{Corollary}
\newtheorem{lemma}[proposition]{Lemma}
\theoremstyle{definition}

\newtheorem{definition}{Definition}[section]
\newtheorem*{theorem*}{Theorem}

\theoremstyle{definition}
\newtheorem{remark}{Remark}[section]

\numberwithin{equation}{section}

\newcommand{\abs}[1]{\left| #1 \right|} 
\newcommand{\norm}[1]{\left\|#1\right\|} 
\newcommand{\normtwo}[1]{
{\left\vert\kern-0.25ex\left\vert\kern-0.25ex\left\vert #1
    \right\vert\kern-0.25ex\right\vert\kern-0.25ex\right\vert} }

\newcommand{\om}{\omega}

\makeatletter
\newsavebox\myboxA
\newsavebox\myboxB
\newlength\mylenA

\newcommand*\xoverline[2][0.75]{%
    \sbox{\myboxA}{$\m@th#2$}%
    \setbox\myboxB\null
    \ht\myboxB=\ht\myboxA%
    \dp\myboxB=\dp\myboxA%
    \wd\myboxB=#1\wd\myboxA
    \sbox\myboxB{$\m@th\overline{\copy\myboxB}$}
    \setlength\mylenA{\the\wd\myboxA}
    \addtolength\mylenA{-\the\wd\myboxB}%
    \ifdim\wd\myboxB<\wd\myboxA%
       \rlap{\hskip 0.5\mylenA\usebox\myboxB}{\usebox\myboxA}%
    \else
        \hskip -0.5\mylenA\rlap{\usebox\myboxA}{\hskip 0.5\mylenA\usebox\myboxB}%
    \fi}
\makeatother

\newcommand{\Proj}{\mathbb{P}(\R^m)}

\newcommand{\Gr}{{\rm Gr}}

\newcommand{\ind}{\mathds{1}}

\newcommand\restr[2]{{
  \left.\kern-\nulldelimiterspace 
  #1 
  \vphantom{\big|} 
  \right|_{#2} 
  }}

\newcommand{\Prob}{\mathrm{Prob}}

\newcommand{\Lip}{\mathrm{Lip}}

\title[Lyapunov exponents of Markov cocycles]{H\"older continuity of the Lyapunov exponent  for Markov cocycles \\via Furstenberg's Formula}

\date{}

\begin{document}

\author[A. Cai]{Ao Cai}
\address{Departamento de Matem\'atica, Pontif\'icia Universidade Cat\'olica do Rio de Janeiro (PUC-Rio), Brazil}
\email{godcaiao@gmail.com}

\author[M. Dur\~aes]{Marcelo Dur\~aes}
\address{Departamento de Matem\'atica, Pontif\'icia Universidade Cat\'olica do Rio de Janeiro (PUC-Rio), Brazil}
\email{accp95@gmail.com}

\author[S. Klein]{Silvius Klein}
\address{Departamento de Matem\'atica, Pontif\'icia Universidade Cat\'olica do Rio de Janeiro (PUC-Rio), Brazil}
\email{silviusk@puc-rio.br}

\author[A. Melo ]{Aline Melo}
\address{Departamento de Matem\'atica, Pontif\'icia Universidade Cat\'olica do Rio de Janeiro (PUC-Rio), Brazil}
\email{alinedemelo.m@gmail.com}

\begin{abstract}
This paper is concerned with the study of linear cocycles over uniformly ergodic Markov shifts on a compact space of symbols. We establish the joint H\"older continuity of the maximal Lyapunov exponent as a function of the cocycle and the transition kernel in the vicinity of any irreducible cocycle with simple maximal Lyapunov exponent. Our approach, via Furstenberg's formula, shows the H\"older continuous dependence on the data of the stationary measure of the projective cocycle and in particular provides a more computable H\"older exponent.   
\end{abstract}

\maketitle


\section{Introduction and statements}\label{intro}
Let $\Sigma$ be a compact metric space. A Markov transition kernel on $\Sigma$ is any continuous map $K \colon \Sigma \to \Prob (\Sigma)$, where the set $\Prob (\Sigma)$ of probabilities on $\Sigma$ is endowed with the weak* topology. 

The iterated Markov kernels $K^n$, $n\ge 1$ are defined inductively by $K^1 = K$ and $K^{n+1}_x (E) = \int_\Sigma K^n_y (E) \, d K_x (y)$, for all $x \in \Sigma$ and all Borel sets $E \subset \Sigma$. We assume that the kernel $K$ is uniformly ergodic, in the sense that for some $n_0 \in \N$ and $\sigma \in (0, 1)$, the total variation norm $\norm{K_x^{n_0}}_{{\rm TV}} \le \sigma$ for all $x\in \Sigma$. This in particular implies the uniqueness (the existence is guaranteed by general principles) of the $K$-stationary measure, that is, a measure $\mu \in \Prob (\Sigma)$ such that $\mu (E) = \int_\Sigma K_x (E) \, d \mu (x)$ for all Borel sets $E$. Uniform ergodicity of the kernel $K$ is then equivalent to the exponential and uniform (in $x \in \Sigma$) convergence of $K_x^n$ to $\mu$ relative to the total variation distance. We will refer to the pair $(K, \mu)$ as a Markov system.

Let $\Pp = \Pp_K = \Pp _{(K, \mu)}$ denote the Markov measure on $X^+ = \Sigma^\N$ with initial distribution $\mu$ and transition kernel $K$. We use the same notation for its extension to the space $X = \Sigma^\Z$ of double sided sequences. 
Let $\sigma$ be the forward shift on $X^+$ and on $X$. Then $(X^+, \Pp, \sigma)$ is a measure preserving (non invertible) dynamical system,  
$(X, \Pp, \sigma)$ is its natural invertible extension which we call a Markov shift. It generalizes the sub-shift of finite type (given by a primitive transition matrix).  

A measurable function $A \colon \Sigma \times \Sigma \to \GL_m (\R)$ induces the skew-product dynamical system $F = F_{(A, K)} \colon X \times \R^m \to X \times \R^m$, 
$$F (\om, v) = \left(\sigma \om, A (\om_1, \om_0) v \right) .$$ 

That is, $F_{(A, K)}$ is a linear cocycle over the base dynamics $(X, \Pp _{(K, \mu)}, \sigma)$, where the fiber dynamics is induced by the map $A$. 
We refer to such a dynamical system as a Markov cocycle.

For simplicity we identify the Markov cocycle $F = F_{(A, K)}$ with the pair $(A, K)$. Its iterates are given by
$$F^n (\om, v) = \left(\sigma^n \om, A^n (\om)  v \right) ,$$ 
where for $\om = \{\om_n\}_{n\in\Z} \in X$, 
$$A^n (\om) = A (\om_n, \om_{n-1}) \cdots A (\om_2, \om_1) \, A (\om_1, \om_0) \, .$$

By Kingman's ergodic theorem, the geometric averages of the fiber iterates of the cocycle $F_{(A, K)}$ converge $\Pp _{(K, \mu)}$-a.s.
$$\frac{1}{n} \, \log \, \norm{A^n (\om)} \to L_1 (A, K) $$
and the limit $L_1 (A, K)$ is called the maximal Lyapunov exponent of the system. Replacing the norm (or largest singular value) of the iterates $A^n (\om)$ by the other singular values, we obtain all the other Lyapunov exponents $L_2 (A, K), \ldots, L_m (A, K)$ of the cocycle $(A, K)$.

\medskip

An important problem in ergodic theory concerns the regularity of the Lyapunov exponents as functions of the input data. It turns out that the type of base dynamics and the topology of the space of cocycles greatly influence the kind of regularity of the Lyapunov exponents, or lack thereof, see~\cite{Bochi-ETDS, Viana-book, DK-book, WangYou} for a sample of available results. 

\medskip

The main goal of this paper is to study the continuity of the maximal Lyapunov exponent of Markov cocycles $F_{(A, K)}$ as a function of the fiber map $A$ and the transition kernel $K$.
Consider the set of Markov cocycles
\begin{align*}
\cocycles := \{ (A,K)\colon &A \colon \Sigma \times \Sigma \to \GL_m (\R) \text{ is Lipschitz continuous and } \\
 &K \colon \Sigma \to \Prob(\Sigma) \text{ is uniformly ergodic and}\\ &\text{continuous in the weak* topology}. \}
\end{align*}

This set is naturally endowed with a metric as follows:
$$
d( (A,K), (B,L) ) := \max \{ d_\infty(A,B), d_{W_1}(K,L) \} ,
$$
where if $A, B \in \Lip(\Sigma \times \Sigma , \GL_m(\R))$ are two Lipschitz continuous fiber maps, 
$$
d_\infty(A,B) := \sup_{\omega_0, \omega_1 \in \Sigma} \norm{ A(\omega_0, \omega_1) - B(\omega_0, \omega_1) }
$$
and the distance between two Markov kernels $K, L$ is defined as
$$
d_{W_1}(K, L) = \sup_{w_0 \in \Sigma} W_1(K_{w_0}, L_{w_0}) ,
$$
where $W_1$ is the Wasserstein distance  in the space of probability measures $\Prob(\Sigma )$. Recall that this distance metrizes the weak* topology.

We assume that the cocycle  $(A, K)$ is quasi-irreducible with respect to the Markov system $(K, \mu)$ (which is a generic condition). Irreducibility refers to the non-existence of a proper, $A$-invariant section, that is, of a measurable function  $V \colon \Sigma \to \Gr(\R^{m})$ (here  $\Gr(\R^{m})$ denotes the Grassmannian of $\R^m$) such that $0 < \dim (V) < m$ and
$$
A(\omega_{n+1},\omega_{n}) V (\omega_{n}) = V (\omega_{n+1}), \; \text{for} \; \mathbb{P}_{(K, \mu)}\text{-a.e.} \; \omega=\{\omega_n\}_n \, .
$$
Quasi-irreducibility is a weaker version of this property, where such a proper $A$-invariant section $V$ may be allowed to exist, but in this case, the maximal Lyapunov exponent of the fiber restriction of the cocycle $F_{(A, K)}$ along the invariant section $V (\om_0)$ must equal $L_1 (A, K)$ for all $\om_0 \in \Sigma$.

We are now ready to formulate the main result of this paper.

\begin{theorem} \label{main thm intro}
Let $(A,K) \in \mathcal{C}$ and assume that:
\begin{enumerate}[(i)]
\item $A$ is quasi irreducible with respect to $(K,\mu)$,
\item $L_1(A,K) > L_2(A,K)$.
\end{enumerate}
Then there exists a neighborhood of $(A,K)$ in $(\mathcal{C},d)$ where the map $(B, L) \mapsto L_1 (B, L)$ is H\"older continuous.
\end{theorem}

\begin{remark}
Assume that all exterior powers $\wedge_k A$, $1\le k \le m$ of the cocycle $A$ are quasi irreducible. From this result we then derive the H\"older continuity of all the other Lyapunov exponents, as long as they are simple. This in particular implies the continuity (not necessarily H\"older) of all exponents, regardless of them being simple or not.
\end{remark}

This result extends~\cite[Theorem 5.1]{DK-book}, where it was established the H\"older continuity of the Lyapunov exponents with respect to the fiber map $A$. In the present work we also allow the transition kernel $K$ to vary, and prove the joint H\"older continuity in $(A, K)$ of the exponents. Moreover, the approach used in this paper (which we briefly explain below) is different from the one in~\cite{DK-book} (which first establishes uniform large deviations type estimates for the cocycle and then deduces the H\"older continuity of the exponents from an abstract continuity theorem). The advantage of the method employed here, besides being more straightforward, is that it provides a more explicit, computable, value of the H\"older exponent (see Remark~\ref{compute}). 

In the case when the space $\Sigma$ of symbols is finite, that is, when the base dynamics is a sub-shift of finite type, there are other results available. 
Fixing any fiber map, the maximal Lyapunov exponent depends analytically on the transition probabilities, see~\cite{Peres2}. This suggests that in our more general setting, the regularity with respect to the transition kernel $K$ might be much higher.  
Moreover, the continuity of the Lyapunov exponent (but without a modulus of continuity) was established in~\cite{Viana-M} without any irreducibility assumption, assuming that the fiber dynamics is two dimensional and depends on only one coordinate.

All of these results, including the one in this paper, are in part inspired by the seminal works of Furstenberg and Kifer~\cite{Furstenberg-K}, Le Page~\cite{LePage} and Peres~\cite{Per}.

\medskip

Let us describe our approach for proving Theorem~\ref{main thm intro}. Let $\Proj$ denote the projective space of dimension $m$ over $\R$ and for $v \in \R^m$, $v\neq 0$, let $\hat v \in \Proj$ be the corresponding projective point. Given a Markov linear cocycle $F = (A, K) \colon X \times \R^m  \to X \times \R^m$, denote by $\hat F = (\hat A, K) \colon X \times \Proj  \to X \times \Proj$ its projective counterpart, which determines the Markov chain on $\Sigma \times \Sigma \times \Proj$,
\begin{equation}\label{chain intro}
\left(\om_1, \om_0, \hat v \right) \to \left(\om_2, \om_1, \hat A (\om_1, \om_0) \hat v \right) \to \cdots
\end{equation}

It turns out that under the assumptions of Theorem~\ref{main thm intro}, the transition kernel $\bar K$ of this Markov chain has a unique stationary measure $m_{(A, K)}$. By Furstenberg's formula, it then follows that the maximal Lyapunov exponent of the cocycle $(A, K)$ is given by
\begin{equation*}
L_1 (A, K) = \int_{ \Sigma \times \Sigma \times \Proj } \, \log \norm{A (\om_1, \om_0) \,  \frac{v}{\norm v}} \, d m_{(A, K)} (\om_1, \om_0, \hat v) \, .
\end{equation*}

Establishing continuity properties of the maximal Lyapunov exponent can thus be reduced to understanding the dependence of the stationary measure $m_{(A, K)}$ on the cocycle $(A, K)$.

Consider the Markov operator $\bar{Q} = \bar{Q}_{(A, K)}$ on $C^0(\Sigma \times \Sigma \times \Proj)$ corresponding to the transition kernel $\bar K$:
 $$
 (\bar{Q}_{(A, K)} \psi)(\omega_1,\omega_0, \hat{v})=\int \psi(\omega_2, \omega_1, \hat{A}(\omega_1,\omega_0)\hat{v}) \; dK_{\omega_1}(\omega_2) .
 $$
 Under the assumptions of Theorem~\ref{main thm intro}, the Markov operator $\bar{Q}_{(A, K)}$ will be shown to be quasi-compact and simple on an appropriate space of observables (that contains the function appearing in Furstenberg's formula). This is equivalent to showing the following convergence (in an appropriate sense) of the powers of the Markov operator:
 $$\bar{Q}_{(A, K)}^n \psi \to \int \psi \, d m_{(A, K)} \quad \text{as } n \to \infty$$
for all such observables $\psi$, which we refer to as strong mixing.

Moreover, for each $n$, the map
$(A, K) \mapsto \bar{Q}_{(A, K)}^n \psi$
will turn out to be H\"older continuous with the same H\"older parameters for all $n$, thus ensuring the same property for the limit quantity as $n\to \infty$.
Choosing $\psi (\om_1, \om_0, \hat v) := \log \norm{A (\om_1, \om_0) \,  \frac{v}{\norm v}}$, we obtain the H\"older continuity of the Lyapunov exponent via Furstenberg's formula.

We note that the strong mixing property of the Markov operator also implies statistical properties such as large deviations and a central limit theorem for the Markov operator. This is due to a recent abstract large deviations estimate result for Markov processes and a  central limit theorem of Gordin and Liv\v{s}ic (see~\cite{CDK-paper3}). These results are not new, they were previously obtained, by different methods, in~\cite{DK-book} and respectively~\cite{Bou}.

\medskip

The rest of the paper is organized as follows. In Section~\ref{preliminar} we formally introduce the Markov operator, the stationary measure and study their basic properties.
In Section~\ref{kifer} we establish the Kifer non-random filtration for Markov cocycles, a more precise version of Oseledets theorem in this context. As a corollary, under the assumption of quasi-irreducibility,  we obtain the uniform convergence of the expected value of the finite scale directional Lyapunov exponent. This is then used in Section~\ref{unique} to establish the strong mixing of the Markov operator and the convergence (in an appropriate sense) of its powers to the unique stationary measure. Finally, in Section~\ref{joint} we obtain the H\"older continuity of the Lyapunov exponent via Furstenberg's Formula.

\section{Markov operators and stationary measures}\label{preliminar}

Let $\Sigma$ be a compact metric space and let $(K, \mu)$ be a Markov system. That is, $K\colon \Sigma \to \Prob(\Sigma)$, $\omega_0\mapsto K_{\omega_0}$ is continuous and uniformly ergodic, while $\mu \in \Prob (\Sigma)$ is its unique stationary measure in the sense that $\mu = K\ast \mu$. More precisely,
$$
\mu(E)=\int_\Sigma K_x(E)d\mu(x) \quad \forall\, E\subset \Sigma \quad\mbox{measurable}.
$$

\begin{remark}\label{equivalence}Recall that the iterates of a transition kernel $K$ are defined recursively setting $K^1 := K$ and for $n \geq 2, E \in \mathcal{F}$,
$K^n_x(E) := \int_X K^{n-1}_y(E) dK_x(y)$.
The uniform ergodicity of  $K$ is equivalent to the convergence $K^n_{\omega_0} \to \mu$ uniformly in $\omega_0 \in \Sigma$. which is equivalent to the existence of constants $\sigma \in (0, 1)$ and $C<\infty$ such that for all $\omega_0\in \Sigma$
	$$
	\norm{K_{\omega_0}^n - \mu}_{{\rm TV}} \leq C\sigma^n \, .
	$$
Furthermore, this is equivalent to the following: $\forall \,\varphi \in L^\infty(\mu)$, 
	$$
	\norm{Q^n\varphi-\int \varphi d\mu}_\infty \leq C\sigma^n \norm{\varphi}_\infty
	$$
	where $Q$ is the Markov operator associated with $K$.
	For these and other characterizations  of uniform ergodicity see~\cite[Theorem 16.0.2]{meyn_tweedie_glynn_2009}.
\end{remark}

Let $A\colon \Sigma \times \Sigma \to \GL_d(\R)$ be a fiber map, which together with the kernel $K$ defines the Markov cocycle $(A, K)$.
We associate to this linear cocycle the transition kernel $\bar{K}\colon \Sigma \times \Sigma \times \Proj\to \Prob(\Sigma \times \Sigma \times \Proj)$ given by
 \begin{equation} \label{barK}
 \bar{K}(\omega_1,\omega_0,\hat{v})=K_{\omega_1}\times \delta_{(\omega_1, \hat{A}(\omega_1,\omega_0)\hat{v})} \, .
 \end{equation}
 The corresponding Markov operator $\bar{Q}$ is defined, for every $\psi \in C^0(\Sigma \times \Sigma \times \Proj)$, by
 $$
 (\bar{Q}\psi)(\omega_1,\omega_0, \hat{v})=\int \psi(\omega_2, \omega_1, \hat{A}(\omega_1,\omega_0)\hat{v}) \; dK_{\omega_1}(\omega_2).
 $$

 Similarly, define the Markov kernel $K_A\colon \Sigma \times \Proj \to \Prob(\Sigma \times \Proj)$,
 \begin{equation} \label{KernelKA}
 K_A(\omega_0,\hat{v})(\cdot)=K_{\omega_0}(\cdot)\times \delta_{\hat{A}(\cdot, \omega_0)\hat{v}}
 \end{equation}
and consider the corresponding Markov operator $Q$, defined, for every $\phi \in C^0(\Sigma \times \Proj)$ by
 $$
 (Q\phi)(\omega_0, \hat{v})=\int \phi(\omega_1, \hat{A}(\omega_1,\omega_0)\hat{v}) \; dK_{\omega_0}(\omega_1).
 $$
 
Moreover, we consider the projection $\Pi\colon C^0(\Sigma \times \Sigma \times \Proj) \to C^0(\Sigma \times \Proj)$:
$$
\Pi\psi(\omega_0,\hat{v})=\int \psi(\omega_1,\omega_0,\hat{v}) \; dK_{\omega_0}(\omega_1).
$$

The following lemma relates these two Markov operators.

\begin{lemma}
	With the notations above, we have $\Pi \circ \bar{Q}=Q\circ \Pi$. That is, the following diagram is commutative.
\[ \begin{tikzcd}
C^0(\Sigma \times \Sigma \times \Proj) \arrow{r}{\bar{Q}} \arrow[swap]{d}{\Pi} & C^0(\Sigma \times \Sigma \times \Proj) \arrow{d}{\Pi} \\%
C^0(\Sigma \times \Proj) \arrow{r}{Q} & C^0(\Sigma \times \Proj)
\end{tikzcd}
\]
\end{lemma}

\begin{proof}
	Take any $\psi \in C^0(\Sigma \times \Sigma \times \Proj)$. A direct computation shows that
	\begin{align*}
		Q\circ \Pi \psi(\omega_0,\hat{v}) &= \int \Pi\psi(\omega_1, \hat{A}(\omega_1,\omega_0)\hat{v}) \; dK_{\omega_0}(\omega_1)\\
		&=\int \int \psi(\omega_2,\omega_1, \hat{A}(\omega_1,\omega_0)\hat{v}) \; dK_{\omega_1}(\omega_2) \, dK_{\omega_0}(\omega_1).
	\end{align*}
On the other hand,
	\begin{align*}
		\Pi\circ \bar{Q}\psi (\omega_0,\hat{v})&=\int \bar{Q}\psi(\omega_1,\omega_0,\hat{v}) \; dK_{\omega_0}(\omega_1)\\
		&=\int \int\psi(\omega_2,\omega_1, \hat{A}(\omega_1,\omega_0)\hat{v}) \; dK_{\omega_1}(\omega_2) \, dK_{\omega_0}(\omega_1).
\end{align*}
This shows $\Pi \circ \bar{Q}=Q\circ \Pi$.
\end{proof}

We denote by $\Prob_{\bar{Q}}(\Sigma\times\Sigma\times\Proj)$ and $\Prob_Q(\Sigma\times\Proj)$, respectively, the convex and compact (since $\Sigma$ is compact) subspace of all $\bar{K}$-stationary probability measure on $\Sigma\times\Sigma\times\Proj$ and $K_A$-stationary probability measure on $\Sigma\times\Proj$.

The following proposition ensures the existence of a $Q$-stationary probability measure in $\Prob(\Sigma\times\Proj)$, given a $\bar{Q}$-stationary probability measure in $\Prob(\Sigma\times\Sigma\times\Proj)$.
 
\begin{proposition} \label{LemmaProductMeasure}
Given $m \in \Prob_{\bar{Q}}(\Sigma\times\Sigma\times\Proj)$, there exists a unique $\eta \in \Prob(\Sigma\times\Proj)$ such that $m = K \ltimes \eta$, that is, for every $\psi \in C(\Sigma \times \Sigma \times \Proj)$,
 $$
 \int \psi(\omega_1,\omega_0, \hat{v}) \, dm(\omega_1,\omega_0,\hat{v}) = \int \psi(\omega_1,\omega_0,\hat{v}) \, dK_{\omega_0}(\omega_1) \, d\eta(\omega_0, \hat{v}).
 $$
Moreover, $\eta$ is a $Q$-stationary probability measure.
\end{proposition}

\begin{proof}
Let $m$ be a $\bar{Q}$-stationary probability measure on $\Sigma\times\Sigma\times\Proj$. Note that if $\psi_1,\psi_2 \in C(\Sigma\times \Sigma \times \Proj)$ are such that $\Pi \psi_1 = \Pi \psi_2$ then $\int \psi_1 \, dm = \int \psi_2 \, dm$.

By Riesz-Markov-Katutani's Theorem, there exists a unique probability measure $\eta$ in $\Prob(\Sigma\times \Proj)$ such that for every $\psi \in C(\Sigma \times \Sigma \times \Proj)$,
	$$
	\int \psi \; dm = \int \phi \; d\eta, \qquad \text{where} \; \phi = \Pi\circ\psi.
	$$
In other words,
	$$
	\int \psi(\omega_1,\omega_0, \hat{v}) \, dm = \int \psi(\omega_1,\omega_0,\hat{v}) \, dK_{\omega_0}(\omega_1) d\eta(\omega_0, \hat{v}).
	$$
	
Moreover, let $\psi \in C(\Sigma \times \Sigma \times \Proj)$ and $\phi \in C(\Sigma \times \Proj)$ such that $\phi = \Pi \circ \psi$. By the definition of $\eta$ and since $m$ is a $\bar{Q}$-stationary probability measure on $\Sigma \times \Sigma \times \Proj$, note that
	$$
	\langle \phi, \eta \rangle = \langle \phi, Q^*\eta \rangle
	$$
then $\eta$ is a $Q$-stationary probability measure on $\Sigma \times \Proj$.
\end{proof}

\begin{remark}
It turns out that the converse of Proposition \ref{LemmaProductMeasure} is also true, so there is a one-to-one correspondence between $\bar K$-stationary measures and $K_A$-stationary measures. Thus certain properties of the Markov operator $Q$ can easily be transferred to $\bar{Q}$.
\end{remark}

Given $x \in \Sigma$, we denote by $\Pp_x = \Pp_{(K, \delta_x)}$ the Markov measure on $\Sigma^\N$ with initial distribution $\delta_x$ and transition kernel $K$.

Define the projective cocycle $\hat{F}\colon \Sigma^{\N}\times \Proj \to \Sigma^{\N}\times \Proj$ by
	$$
	\hat{F}(\omega, \hat{v}) = (\sigma\omega, \hat{A}(\omega_1, \omega_0)\hat{v}).
	$$
	
\begin{proposition} \label{propexistenciaeta}
Given $\eta \in \Prob_Q(\Sigma \times \Proj)$, there exists an $\hat{F}$-invariant probability measure $\tilde{\eta}$ in $\Prob(\Sigma^{\N}\times \Proj)$ such that
	$$
	(\pi_{01})_* \tilde{\eta} = \eta, \quad \text{where} \quad \pi_{01}(\omega, \hat{v}) = (\omega_1, \omega_0, \hat{v}).
	$$
\end{proposition}

\begin{proof}
Let $\psi \in C(\Sigma^{\N} \times \Proj)$ and define $\tilde{\eta} \in \Prob(\Sigma^{\N}\times \Proj)$ such that
	\begin{equation} \label{etatilde}
	\int \psi(\omega,\hat{v}) \; d\tilde{\eta}(\omega,\hat{v}) := \int \psi(\omega,\hat{v}) \, d\Pp_{\omega_0}(\omega) \, d\eta(\omega_0,\hat{v})
	\end{equation}

A simple calculation shows that $(\pi_{01})_* \tilde{\eta} = \eta$. Moreover, $\tilde{\eta}$ is an $\hat{F}$-invariant measure. In fact, since $\eta$ is a $Q$-stationary probability measure, for every $\phi \in C(\Sigma \times \Proj)$,
	\begin{align*}
		\int \phi(\omega_0,\hat{v})\; d\eta(\omega_0,\hat{v}) 
		&= \int \phi(\omega_1,\hat{A}(\omega_1,\omega_0)\hat{v})\; dK_{\omega_0}(\omega_1) \; d\eta(\omega_0,\hat{v}).
	\end{align*}
Hence, given $\psi \in C(\Sigma^{\N}\times \Proj)$
	\begin{align*}
		\int \psi(\omega,\hat{v})\; d\tilde{\eta}(\omega, \hat{v}) 
		&= \int \psi(\omega,\hat{v})\; d\Pp_{\omega_0}(\omega) \; d\eta(\omega_0,\hat{v}) \\
		&= \int Q\left(\int \psi(\omega,\hat{v})\; d\Pp_{\omega_0}(\omega) \right) \; d\eta(\omega_0,\hat{v}) \\
		&= \int (\psi\circ \hat{F})(\omega,\hat{v})\; d\tilde{\eta}(\omega,\hat{v})
	\end{align*}
proving that $\tilde{\eta}$ is an $\hat{F}$-invariant probability measure on $\Sigma^{\N}\times \Proj$.
\end{proof}

\begin{definition}
An observable $\phi \in L^{\infty}(\Sigma \times \Proj)$ is called \emph{$\eta$-stationary} if $Q \phi(x,\hat{v}) = \phi(x,\hat{v})$ for $\eta$-almost every $(x,\hat{v})$.
\end{definition}

\begin{definition}
A Borel set $F \subset \Sigma \times \Proj$ is \emph{$\eta$-stationary} if the indicator function $\mathbbm{1}_F$ is $\eta$-stationary. That is, $F$ is an $\eta$-stationary set if and only if for $\eta$-almost every $(\omega_0,\hat{v})$, we have $(\omega_1, \hat{A}(\omega_1, \omega_0)\hat{v}) \in F$, for $K_{\omega_0}$-almost every $\omega_1 \in \Sigma$.
\end{definition}

The following proposition proves that the probability measure $\eta$ is an extremal point of $\Prob_{Q}(\Sigma \times \Proj)$ if and only if $\tilde{\eta}$ (defined in (\ref{etatilde})) is an $\hat{F}$-ergodic probability measure. 

\begin{proposition} \label{propequiv}
Given $\eta \in \Prob_{Q}(\Sigma \times \Proj)$, the following are equivalent:

	\begin{itemize}
		\item[(i)] $\eta$ is an extremal point of $\Prob_{Q}(\Sigma \times \Proj)$

		\item[(ii)] If $F \subset \Sigma \times \Proj$ is an $\eta$-stationary set then $\eta(F) = 0$ or $\eta(F) = 1$

		\item[(iii)] If $\phi \in L^{\infty}(\Sigma \times \Proj)$ is an $\eta$-stationary function then $\phi$ is a constant function $\eta$-almost everywhere

		\item[(iv)] The system $(\Sigma^{\N} \times \Proj, \hat{F}, \tilde{\eta})$ is ergodic.
	\end{itemize}
\end{proposition}

\begin{proof}
For the proof of the equivalence (ii) $\Leftrightarrow$ (iii) see~\cite[Proposition 5.11]{Viana-book}.
Let us prove that (i) $\Rightarrow$ (ii). 

Let $\eta$ be an extremal point of $\Prob_{Q}(\Sigma \times \Proj)$ and assume, by contradiction, that there exists an $\eta$-stationary subset $F \subset \Sigma \times \Proj$ such that $t := \eta(F) \in (0,1)$. Then, $F^c := (\Sigma \times \Proj) \setminus F$ is also an $\eta$-stationary subset and $\eta(F^c) = 1-t \in (0,1)$.

Let $\eta_F$ and $\eta_{F^c}$ be probability measures on $\Sigma \times \Proj$ such that
$$\eta_F(E) = \frac{\eta(E \cap F)}{\eta(F)} \quad \text{and} \quad \eta_{F^c}(E) = \frac{\eta(E \cap F^c)}{\eta(F^c)}.$$

Observe that $\eta_F \neq \eta_{F^c}$ and $\eta = t\eta_F + (1-t)\eta_{F^c}$. 
Since the indicator function $\mathbbm{1}_F$ is $\eta$-stationary and $\eta$ is a $Q$-stationary probability measure, we have that for every $\phi \in C(\Sigma \times \Proj)$,
	\begin{align*}
		\int_{\Sigma \times \Proj} Q\phi(\omega_0, \hat{v}) \; d\eta_F(\omega_0, \hat{v}) 
		&= \frac{1}{\eta(F)} \int_{F} Q\phi(\omega_0, \hat{v}) \; d\eta(\omega_0, \hat{v}) \\
		&= \frac{1}{\eta(F)} \int_{F} \phi(\omega_0, \hat{v})  \; d\eta(\omega_0, \hat{v}) \\
		&= \int_{\Sigma \times \Proj} \phi \; d\eta_F
	\end{align*}
and then $\eta_F$ is an $Q$-stationary probability measure. Analogously, $\eta_{F^c}$ is also $Q$-stationary and this contradicts the hypothesis that $\eta$ is an extremal point of $\Prob_{Q}(\Sigma \times \Proj)$.

Suppose now that (iii) is true. We will proof that $\tilde{\eta}$ is an $\hat{F}$-ergodic probability measure. Let $\psi \in L^{\infty}(\Sigma^{\N}\times \Proj)$ such that $\psi \circ \hat{F} = \psi, \tilde{\eta}$-almost everywhere and consider $\phi: \Sigma \times \Proj \to \R$ such that
	$$
	\phi(\omega_0,\hat{v}) = \int \psi(\omega, \hat{v}) \; d\Pp_{\omega_0}(\omega).
	$$
	
Since $\psi$ is an $\hat{F}$-invariant function, it is easy to see that $\phi$ is an $\eta$-stationary function and, consequently, $\phi$ is constant $\eta$-almost surely. Then, $\psi$ is constant in $(\omega, \hat{v})$, $\tilde{\eta}$-almost every $(\omega, \hat{v})$ since $\psi = \psi\circ \hat{F}^k$ and $\phi$ is constant for $\tilde{\eta}$-almost everywhere. This proves (iv).

It remains to proof that (iv) implies (i). Assume by contradiction that $\eta$ is not an extremal point of $\Prob_{Q}(\Sigma \times \Proj)$. Then, there exist $t\in (0,1)$ and $\eta_1, \eta_2 \in \Prob_{Q}(\Sigma \times \Proj)$ such that 
	$$
	\eta = t\eta_1 + (1-t)\eta_2.
	$$
Then, $\tilde{\eta} = t\tilde{\eta}_1 + (1-t)\tilde{\eta}_2$ is not and extremal point and, consequently, $\tilde{\eta}$ is not an ergodic measure.
\end{proof}

\section{Kifer non-random filtration}\label{kifer}

 For measure preserving dynamical systems, the Oseledets multiplicative ergodic theorem improves the Furstenberg-Kesten theorem in that it provides exponential rates of convergence of the iterates $A^{(n)} (\omega)$ of the cocycle $A$ along all directions. Later, Kifer improved the Oseledets theorem for random cocycles (i.e. linear cocycles over a Bernoulli shift) by proving the existence of an invariant filtration that does not depend on the base point, thus it is non-random. The main goal of this section is to derive a version of this result  in the context of Markov cocycles. Bougerol~\cite{Bou} obtained a partial related result, by a different method.
Assuming the quasi-irreducibility of the cocycle we derive a Furstenberg-type formula and eventually the uniform convergence of the expected value of the finite scale Lyapunov exponent.

%
%
%

Given $(K,\mu)$ a Markov system, consider the continuous observable $\psi\colon \Sigma \times \Sigma\times \Proj \to \R$ defined by
	\begin{equation}\label{observablepsi}
	\psi(x,y,\hat{v}) = \log \Vert A(y,x)v \Vert
	\end{equation}
where $v$ is any unit vector representing the projective point $\hat{v}$. The observable $\psi$ extends naturally to a function $\tilde{\psi}\colon \Sigma^{\N} \times \Proj \to \R$ such that $\tilde{\psi} = \psi \circ \pi_{01}$.

Consider the continuous linear functional $\alpha\colon \Prob_{Q}(\Sigma \times \Proj) \to \R$ defined by
	$$
	\alpha(\eta) := \int_{\Sigma \times \Sigma \times \Proj} \psi(y,x,\hat{v})\; dK_x(y)d\eta(x,\hat{v})
	$$
and define
	$$
	\beta := \max\{\alpha(\eta): \eta\in \Prob_{Q}(\Sigma \times \Proj)\} 
	$$
	


Recall the following results of H. Furstenberg and Y. Kifer.

\begin{theorem}[Furstenberg-Kifer \cite{FKi} Theorem 1.1] \label{thmProcMarkovGeral}
Let  $\{Z_n\}_{n\geq 0}$ be a $K$-Markov chain in $\Sigma$   and let $f\in C(\Sigma)$. Then with probability one we have
	$$
	\limsup \dfrac{1}{n} \displaystyle\sum_{j=0}^{n-1} f(Z_j) \leq \sup \left\{\int_{\Sigma} f \; d\nu: \nu \in \Prob_K(\Sigma) \right\}.
	$$
\end{theorem}

\begin{theorem}[Furstenberg-Kifer \cite{FKi} Theorem 1.4] \label{thmProcMarkovConst}
Let  $\{Z_n\}_{n\geq 0}$ be a $K$-Markov chain in $\Sigma$, let $f\in C(\Sigma)$ and assume that $\int_{\Sigma} f\; d\nu = \beta$ for every $K$-stationary probability measure $\nu \in \Prob_K(\Sigma)$. Then with probability one we have
		$$
		\lim \frac{1}{n} \sum_{j=0}^{n-1} f(Z_j) = \beta.
		$$
\end{theorem}

We apply these general results to our setting to derive the following.
	
\begin{theorem} \label{ThmFK}
For every $(\omega_0,v) \in \Sigma \times \R^m$,
	\begin{itemize}

		\item[(i)] $\limsup \frac{1}{n}\log \Vert A^n(\omega)v \Vert \leq \beta$, for $\Pp_{\omega}$-almost every $\omega$

		\item[(ii)] If $\alpha$ is constant then
$$\lim \frac{1}{n}\log \Vert A^n(\omega)v \Vert = \beta$$
with full probability.

		\item[(iii)] $\lim \frac{1}{n}\log \Vert A^n(\omega) \Vert = \beta$, for $\Pp_{\omega_0}$-almost every $\omega$. In particular, 
		$L_1 (A, K) = \beta$, which is a version of Furstenberg's formula in this setting.

	\end{itemize}
\end{theorem}

\begin{proof}
Let $M = \Sigma \times \Sigma \times \Proj$, $\bar{K}: M\to \Prob(M)$ be the kernel defined in (\ref{barK}) and the continuous observable $\psi \in C^0(M)$ defined in (\ref{observablepsi}).
	
For each $v \in \R^m\setminus \{0\}$, consider the $\bar{K}$-Markov random process $Z_n^{\hat{v}}\colon \Sigma^{\N} \to M, n\geq 0$, defined by
	$$
	Z_n^{\hat{v}}(\omega) := (\omega_{n+1}, \omega_{n}, \hat{A}^n(\omega)\hat{v}).
	$$


Items (i) and (ii) follow applying Theorem \ref{thmProcMarkovGeral} and Theorem \ref{thmProcMarkovConst} respectively.

Let us now prove (iii). By Furstenberg-Kesten's Theorem, 
	$$
	\lim \frac{1}{n}\log \Vert A^n(\omega) \Vert = L_1(A), \quad \Pp_{\omega_0}-\text{a.e.} \; \omega \in \Sigma^{\N}
	$$
	
Fixing a basis $\{e_1, \ldots, e_d\}$ of $\R^m$, define the matrix norm
	$$
	\Vert g \Vert'= \max_{1\leq i\leq d} \Vert ge_i \Vert.
	$$
The set of maximizing measures
	$$
	\mathcal{M} := \{\eta \in \Prob_Q(\Sigma\times\Proj): \alpha(\eta) = \beta\}
	$$
is a non-empty compact convex set. By Krein-Milman Theorem there exists an extremal point $\eta$ of $\mathcal{M}$ and then this measure is also an extremal point of $\Prob_Q(\Sigma\times\Proj)$. By Proposition \ref{propequiv}, $\tilde{\eta}$ is an $\hat{F}_A$-ergodic probability measure. Thus, by Birkhoff Ergodic Theorem, for $\tilde{\eta}$-almost every $(\omega, \hat{v}) \in \Sigma^{\N}\times \Pp$,
	\begin{align*}
		\beta = \alpha(\eta) 
		&= \int_{\Sigma \times \Sigma \times \Proj} \psi\; d(K \ltimes \eta) = \int_{\Sigma^{\N} \times \Proj} \tilde{\psi}\; d\tilde{\eta} \\
		&= \lim_{n\to \infty} \frac{1}{n}\sum_{j=0}^{n-1} \tilde{\psi}(\hat{F}^j(\omega, \hat{v})) = \lim_{n\to \infty} \frac{1}{n}\log \Vert A^n(\omega) v \Vert \\
		&\leq \lim_{n\to \infty} \frac{1}{n}\log \Vert A^n(\omega) \Vert = \lim_{n\to \infty} \frac{1}{n}\log \Vert A^n(\omega) \Vert' \\
		&= \max_{1\leq i\leq d} \lim_{n\to \infty} \frac{1}{n}\log \Vert A^n(\omega)e_i \Vert \leq \beta
	\end{align*}
This proves (iii).
\end{proof}

\begin{definition}
Let $\Gr(\R^{m})$ denote the Grassmann manifold of the Euclidean space $\R^{m}$. A measurable section $L\colon \Sigma \to \Gr(\R^{m})$ is called $A$-invariant if 
	$$
	A(\omega_{n+1},\omega_{n})L(\omega_{n}) = L(\omega_{n+1}), \; \text{for} \; \mathbb{P}_{\mu}\text{-a.e.} \; \omega=(\omega_n)_{n\in \N}.
	$$
\end{definition}

We are ready to state and prove a version of Kifer's non random filtration theorem for Markov cocycles.

\begin{theorem} \label{KiferNRF}
Let $(\Sigma, K,\mu)$ be a Markov system. Then for every $\omega_0 \in \Sigma$ there exists a filtration $\mathcal{L} = (\mathcal{L}_{0}, \mathcal{L}_{1}, \ldots, \mathcal{L}_{r})$, where $\mathcal{L}_{j}\colon \Sigma \to \Gr(\R^m)$ is a section,  with $0\leq r \leq m$,
	$$
	\{0\} \subsetneq \mathcal{L}_{r}(\omega_0) \subsetneq \ldots \subsetneq \mathcal{L}_{1}(\omega_0) \subsetneq \mathcal{L}_{0}(\omega_0) = \R^m
	$$
and there are numbers
	$$
	\beta = \beta_{0} > \beta_{1} > \cdots > \beta_{r}
	$$ 
such that for every $0\leq j \leq r$

	\begin{itemize}
		\item[(i)] each section $\mathcal{L}_{j}\colon \Sigma \to \Gr(\R^m)$ is $A$-invariant

		\item[(ii)] for every $v \in \mathcal{L}_{j}(\omega_0) \setminus \mathcal{L}_{j+1}(\omega_0)$ and $\Pp_{\omega_0}$-almost every $\omega \in \Sigma^{\N}$,
$$\lim \frac{1}{n} \log \Vert A^{n}(\omega) v \Vert = \beta_{j}$$


		\item[(iii)] $\{\beta_{j}: 0 \leq j \leq r\}$ $= \{\alpha(\eta): \eta \; \text{is an extremal point of} \; \Prob_{Q}(\Sigma \times \Proj)\}$ 

		\item[(iv)] for any extremal point $\eta$ of $\Prob_{Q}(\Sigma \times \Proj)$ such that $\alpha(\eta) = \beta_{j}$, we have $\eta\{(\omega_0, \hat{v}): \hat{v} \in \hat{\mathcal{L}}_j(\omega_0)\} = 1$ and $\eta\{(\omega_0, \hat{v}): \hat{v} \in \hat{\mathcal{L}}_{j+1}(\omega_0)\} = 0$.
	\end{itemize}
\end{theorem}

\begin{proof} 
Let $\omega_0\in \Sigma$ and $\beta= \beta_{0} > \beta_{1} > \cdots > \beta_{r}$ be the finite set $\Sigma(\alpha)$ of all values of the linear functional $\alpha$ over the extremal points of $\Prob_{Q}(\Sigma \times \Proj)$.

\begin{lemma}
The set $\Sigma(\alpha)$ is finite and is contained in the Lyapunov spectrum.
\end{lemma}

\begin{proof} Let $\psi\colon \Sigma \times \Sigma\times \Proj \to \R$ be the continuous observable defined in (\ref{observablepsi}) and its natural extension $\tilde{\psi}\colon \Sigma^{\N} \times \Proj \to \R$, 
$\tilde{\psi} = \psi \circ \pi_{01}$.

By Proposition \ref{propexistenciaeta} and Proposition \ref{propequiv}, there exists an $\hat{F}$-invariant and ergodic probability measure $\tilde{\eta} \in \Prob(\Sigma^{\N} \times \Proj)$ such that $(\pi_{01})_* \tilde{\eta} = \eta$. By Birkhoff Ergodic Theorem, for $\tilde{\eta}$-almost every $(\omega, \hat{v}) \in \Sigma^{\N} \times \Proj$
	\begin{align*}
		\alpha(\eta) 
		&= \int_{\Sigma \times \Sigma \times \Proj} \psi\; d(K\ltimes \eta) = \int_{\Sigma^{\N} \times \Proj} \tilde{\psi} \, d\tilde{\eta} \\
		&= \lim_{n\to \infty} \frac{1}{n}\sum_{j=0}^{n-1} \tilde{\psi}(\hat{F}^j(\omega, \hat{v}))
		= \lim_{n\to \infty} \frac{1}{n} \log \Vert A^n(\omega) v \Vert
	\end{align*}
matches one of the Lyapunov exponents. Hence, $\Sigma(\alpha)$ is a finite set and is contained in the Lyapunov spectrum.
\end{proof}

By Theorem \ref{ThmFK}, $\beta = L_1(A)$. If $\Sigma(\alpha) = \{\beta\}$ then the trivial filtration $\mathcal{L} = \{\mathcal{L}_{0}(\omega_0) = \R^m\}$, with $r=0$, satisfies the conclusions (1)-(4).
Assume now that $r\geq 1$ and consider an extremal point $\eta \in \Prob_{Q}(\Sigma \times \Proj)$ such that $\alpha(\eta) = \beta_{1} < \beta_{0} = \beta$.

\begin{lemma} \label{LemmaLproper}
The set
$$L(\omega_0) := \left\{v\in \R^m: \limsup_{n\to\infty} \frac{1}{n}\log \Vert A^n(\omega)v \Vert \leq \beta_{1} \quad \mathbb{P}_{\omega_0}\text{-a.e} \; \omega\right\}$$
is a proper linear subspace such that $\eta\{(\omega_0,\hat{v}): \hat{v} \in \hat{L}(\omega_0)\} = 1$.
\end{lemma}

\begin{proof}
A direct computation and since $\beta$ is a Lyapunov exponent, we have that $L(\omega_0)$ is a proper linear subspace.

Moreover, since
$$\lim_{n\to\infty} \frac{1}{n} \log \Vert A^n(\omega)v \Vert = \beta_{1} = \alpha(\eta), \quad \tilde{\eta}\text{-a.e.} \; (\omega, \hat{v})$$
we conclude that $\eta\{(\omega_0, \hat{v}): \hat{v} \in \hat{L}(\omega_0)\} = 1$.
\end{proof}


\begin{lemma}
For $\mu$-almost every $x,y \in \Sigma$, $\dim L(x) = \dim L(y)$.
\end{lemma}

\begin{proof}
Define
$$E_v := \left\{\omega\in \Sigma^{\N}: \limsup_{n\to\infty} \frac{1}{n}\log \Vert A^n(\omega)v \Vert \leq \beta_{1} \right\}.$$

Note that $v \in L(x)$ if, and only if, $\Pp_{x}(E_v) = 1$. For every $x, y\in \Sigma$, consider $\Pp_{yx} \in \Prob(\Sigma^{\N})$ such that 
$$\Pp_{x} := \int_{\Sigma} \Pp_{yx} \; dK_x(y).$$

If $v\in L(x)$,
$$1 = \Pp_{x}(E_v) = \int_{\Sigma} \Pp_{yx}(E_v) \; dK_x(y)$$
then $\Pp_{yx}(E_v) = 1$ for $K_x$-almost every $y$. Hence,
\begin{equation*}
	\limsup\frac{1}{n} \log \Vert A(\omega_n, \omega_{n-1})\cdots A(\omega_2, y) A(y,x)v \Vert \leq \beta_1 \quad \Pp_{yx}\text{- a.e.} \; \omega. 
\end{equation*} 
and, consequently, $\bar{v} = A(y,x)v \in L(y)$ for $K_x$-almost every $y$.
This implies that $\dim L(x) \leq \dim L(y)$ for $K_x$-almost every $y$ and then $K_x(E) = 1$, where
$$E := \{y \in \Sigma: \dim L(x) \leq \dim L(y)\}.$$

We can prove by induction that $K^n_x(E) = 1$ for every $n\in \N$ and $x\in \Sigma$. Since $K^n_x \to \mu$ it follows that $\mu(E) = 1$ and this conclude the proof.
\end{proof}

Consider all sections $\mathfrak{L} = \{(x, L(x)): x\in \Sigma\}$ with $\dim L(x)$ constant in $x$ and $\eta(\mathfrak{\hat{L}}) = 1$, where $\mathfrak{\hat{L}} = \{(x, \hat{v}): x\in \Sigma, \hat{v} \in \hat{L}(x)\}$. Let $\mathfrak{L}_{\eta} \subset \mathfrak{L}$ the minimal section, in the sense that $L_{\eta}(x) := L(x)$ is the minimal subspace with respect to the dimension.

\begin{lemma} \label{lemmasectioninvariant}
The section $L_{\eta}\colon \Sigma \to \Gr(\R^m)$ is $A$-invariant.
\end{lemma}

\begin{proof}
Define $E := \{(\omega_0, \hat{v}): \hat{v} \in \hat{L}(\omega_0)\}$. Since $\eta$ is a stationary measure and $\eta\{E\} = 1$, we can show that
	$$
	\int \mathbbm{1}_{E}(\omega_1, \hat{A}(\omega_1, \omega_0)\hat{v}) \; dK_{\omega_0}(\omega_1) = 1.
	$$
That is, for every $\omega_0 \in \Sigma$ and for $K_{\omega_0}$-almost every $\omega_1 \in \Sigma$,
\begin{align*}
	(\omega_1, \hat{A}(\omega_1, \omega_0)\hat{v}) \in E 
	&\Leftrightarrow \hat{A}(\omega_1, \omega_0)\hat{v} \in \hat{L}_{\eta}(\omega_1).
\end{align*}
Hence, the section $L_{\eta}: \Sigma \to \Gr(\R^m)$ is $A$-invariant.
\end{proof}

Given an $A$-invariant section $L\colon \Sigma \to \Gr(\R^m)$, consider the fiber spaces $\mathfrak{L} = \cup_{\omega \in \Sigma^{\N}} (\{\omega\} \times L(\omega_0))$ and $\R^m / \mathfrak{L} := \cup_{\omega \in \Sigma^{\N}} (\{\omega\} \times \R^m / L(\omega_0))$.  Then, it induces the fiber bundles $F_{L} := F|_{\mathfrak{L}}$ and $F_{\R^m/L} := F|_{\R^m /\mathfrak{L}}$. Consider the fiber bundles $\mathfrak{L}_{0} = \cup_{\omega_0 \in \Sigma} (\{\omega_0\} \times L(\omega_0))$ and $(\R^m / \mathfrak{L})_{0} := \cup_{\omega_0 \in \Sigma} (\{\omega_0\} \times \R^m / L(\omega_0))$. 
We denote by $\beta_{L}$ and $\beta_{\R^m / L}$ the maxima of the linear functionals $\alpha: \Prob(\mathfrak{L}_{0}) \to \R$ and $\hat{\alpha}: \Prob((\R^m / \mathfrak{L})_{0}) \to \R$ respectively.

\begin{proposition} \label{betamax}
Let $L\colon \Sigma \to \Gr(\R^m)$ be an $A$-invariant section. Then
$$\beta = \max\{\beta_L, \beta_{\R^m / L}\}.$$
\end{proposition}

\begin{proof}
Without loss generality we can assume that $L(\omega_0) = \R^k \times \{0\}$ and so we can write $A: \Sigma^{\N} \to \GL_d(\R)$ in the form
$$A = \left[\begin{array}{cc}
B & C \\ 
0 & D
\end{array}\right] $$
with block components $B\colon \Sigma^{\N} \to \GL_k(\R), D\colon \Sigma^{\N} \to \GL_{d-k}(\R)$ and $C\colon \Sigma^{\N} \to \Mat_{k \times (d-k)}(\R)$.

The functions $B$ and $D$ represent the fiber bundles associates with $F_L$ and $F_{\R^m/L}$ respectively. A simple calculation gives
$$A^n(\omega) =
\left[
\begin{array}{cc}
B^n(\omega) & C_n(\omega) \\ 
0 & D^n(\omega)
\end{array} 
\right]$$
where
$$C_n := \sum_{j=0}^{n-1} B^{n-i-1}(T^i\omega) C(T^i\omega) D^i(\omega).$$

Since $\max\{\Vert B^n \Vert, \Vert D^n \Vert\} \leq \Vert A^n \Vert$ we have
$$\max\left\{\limsup \frac{1}{n}\log\Vert B^n \Vert, \limsup \frac{1}{n}\log\Vert D^n \Vert\right\} \leq \limsup \frac{1}{n}\log\Vert A^n \Vert$$
which implies that $\max\{\beta_L, \beta_{\R^m/L}\} \leq \beta$. On the other hand, since $C$ is bounded, the above formula shows that $\Vert A^n \Vert$ can never grow exponentially faster than both $\Vert B^n \Vert$ and $\Vert D^n \Vert$. Thus,
$$\limsup \frac{1}{n}\log\Vert A^n \Vert \leq \max\left\{\limsup \frac{1}{n}\log\Vert B^n \Vert, \limsup \frac{1}{n}\log\Vert D^n \Vert\right\}$$
and then $\beta \leq \max\{\beta_L, \beta_{\R^m/L}\}$.
\end{proof}

Observe that $L_{\eta}\colon \Sigma \to \Gr(\R^m)$ is a section with the following properties:

\begin{itemize}

\item[(a)] $L_{\eta}$ is an $A$-invariant section,

\item[(b)] $L_{\eta}$ has constant dimension of each fiber,

\item[(c)] $\eta\{(\omega_0, \hat{v}): \hat{v} \in L_{\eta}(\omega_0)\} = 1$.

\end{itemize}

Amongst all sections $L\colon \Sigma \to \Gr(\R^m)$ such that (a), (b) and (c) hold, choose a section $L_1\colon \Sigma \to \Gr(\R^m)$ that has maximal dimension, $L_{\eta} \subset L_1$ and $\beta_{L_{1}} < \beta$.

Note that $L_1$ is an $A$-invariant section such that $\beta_{L_1} = \beta_{1} < \beta$ and $L_1(\omega_0)$ is a proper linear subspace of $\R^m$. Then, by Proposition \ref{betamax}, we must have $\beta = \beta_{\R^m/L_1}$. 

So far we have proven (i), (iii) and (iv) for the section $L_1$ and the constant section $L_0 = \R^m$. Next goal is to prove (ii) for vectors in $L_0(\omega_0) \setminus L_1(\omega_0)$.

\begin{proposition} \label{betaconst}
For any $\hat{\theta}$ stationary measure in $\Prob((\R^m/\mathfrak{L})_0)$, we have $\hat{\alpha}(\theta) = \beta_{\R^m/L_1}$.
\end{proposition}

\begin{proof}
Suppose that there exists a stationary measure $\hat{\theta}$ in $\Prob((\R^m/\mathfrak{L})_0)$ such that $\hat{\alpha}(\theta) \neq \beta_{\R^m/L_1}$. Since $\beta_{\R^m/L_1}$ is the maxima of the linear functional $\hat{\alpha}\colon \Prob((\R^m / M)_{0}) \to \R$, we have that $\hat{\alpha}(\theta) < \beta_{\R^m/L_1} = \beta$. Hence, as we have seen above, there exists an $A$-invariant proper section $L'\colon \Sigma \to \Gr(\R^m)$ such that $L_1(\omega_0) \subsetneq L'(\omega_0) \subsetneq \R^m$ and $\beta_{L'} < \beta$. This contradicts the maximality of the section $L_1$.
\end{proof}

\begin{proposition}
For each $v \in \R^m \setminus L_1(\omega_0)$ and $\Pp_{\omega_0}$-almost every $\omega \in \Sigma^{\N}$,
$$\lim \frac{1}{n} \log \Vert A^n(\omega)v \Vert = \beta.$$
\end{proposition}

\begin{proof}
Since $\beta_{L_1} < \beta$, by Proposition \ref{betamax} and Proposition \ref{betaconst}, the linear functional $\hat{\alpha}\colon \Prob((\R^m / \mathfrak{L})_{0}) \to \R$ is constant. Hence, by Theorem \ref{ThmFK}, for $\Pp_{\omega_0}$-almost every $\omega$
$$\beta = \beta_{\R^m/L_1} = \lim \frac{1}{n}\log\Vert D^n(\omega)v \Vert$$
where $D\colon \Sigma^{\N} \to \GL_{d-k}(\R)$ is the block component of $A$ that represents the quotient cocycle induced by $A$ on $\R^m/\mathcal{L}_1$. Since for any $(0,v) \in \R^m \setminus L_1(\omega_0)$,
$$\Vert A^n(\omega)(0,v) \Vert \geq \Vert D^n(\omega)v \Vert.$$
$$\lim \frac{1}{n} \log \Vert A^n(\omega)v \Vert \geq \lim \frac{1}{n}\log\Vert D^n(\omega)v \Vert = \beta$$
the conclusion follows immediately.
\end{proof}

The proposition above shows that item (ii) holds for the section $L_1\colon \Sigma \to \Gr(\R^m)$. If $r = 1$, by Theorem \ref{ThmFK}, item (ii) holds also for the subspace $L_1(\omega_0)$. Otherwise $r \geq 2$ and we can apply the same procedure to the cocycle induced by $A$ on $\mathfrak{L}_1$ to get another $A$-invariant section $L_2\colon \Sigma \to \Gr(\R^m)$ with $L_2(\omega_0) \subsetneq L_1(\omega_0)$ such that for every $v \in L_1(\omega_0) \setminus L_2(\omega_0)$ and $\Pp$-almost every $\omega$
$$\lim \frac{1}{n} \log \Vert A^n(\omega)v \Vert = \beta_1.$$
The completion of the proof follows by induction.
\end{proof}

Let us introduce some generic notions.
\begin{definition}
A Markov cocycle $A$ is called \emph{quasi-irreducible} w.r.t.$\,(K,\mu)$ if there exists no measurable proper $A$-invariant section $V\colon\Sigma \to \Gr(\R^m)$ such that $L|_{V(\omega_0)}<L_1(A), \forall\, \omega_0\in \Sigma$.
\end{definition}

An immediate consequence of the quasi-irreducibility condition is the following.
\begin{corollary} \label{QITrivialFiltr}
	Assume that $A$ is quasi-irreducible, then the non-random filtration is trivial. Namely, we have 
	$$\mathcal{L}=\{L_0\},\quad\R^m=L_0(\omega_0), \,\forall\, \omega_0 \in \Sigma.$$
	In particular, $\alpha(\eta)\equiv \beta=L_1(A, K)$ is a constant for all $\eta \in \Prob_{Q}(\Sigma\times \Pp(\R^m))$. 
\end{corollary}

\begin{proof}
	Assume by contradiction that the non-random filtration is not trivial. Then $\mathcal{L}$ must contain at least one more element $L_1\neq L_0$ such that $\forall\, \omega_0\in \Sigma$, for every $v\in L_1(\omega_0)$ and $\Pp_{\omega_0}$-a.e. $\omega\in \Sigma^\N$,
	\begin{equation} \label{QITrivialFiltreq}
	\lim_{n\to\infty}\frac{1}{n}\log \norm{A^n(\omega)v} \leq \beta_1<\beta_{0}=L_1(A).
	\end{equation}

	This contradicts to $A$ being quasi-irreducible. This proves the non-random filtration is trivial. By item (iii) of Theorem \ref{KiferNRF}, we have 
	$$
	\alpha(\eta)\equiv \beta_0=\beta=L_1(A)=\lim_{n \to \infty}\frac{1}{n}\log \norm{A^n(\omega)}, \Pp_{\omega_0}\text{-a.e.} \,\omega \in \Sigma^\N
	$$
	for all $\eta$ that is an extremal point of $\Prob_{Q}(\Sigma\times \Pp(\R^m))$. 
	By Krein-Milman and the linearity of $\alpha(\cdot)$, the result follows.
\end{proof}

The following theorem ensures the uniform convergence of the expected value in $\omega_0$.

\begin{theorem} \label{convUnif}
Let $A$ be a Markov cocycle over a Markov system $(K,\mu)$ such that $A$ and $A^{-1}$ are both measurable. If  $A$ is quasi-irreducible and $L_1(A, K) > L_2(A, K)$ then
$$\lim_{n\to\infty} \frac{1}{n} \mathbb{E}_{\omega_0} (\log \Vert A^{(n)}(\omega)v \Vert) = L_1(A,K),$$
with uniform convergence in $(\omega_0,\hat{v})\in \Sigma\times \mathbb{S}^{d-1}$.
\end{theorem}

\begin{proof}
Since $A$ is quasi-irreducible, by Corollary \ref{QITrivialFiltr} and Lebesgue dominated convergence theorem, we have the pointwise convergence:
	\begin{equation} \label{pointconver}
	\lim_{n\to\infty} \mathbb{E}_{\omega_0}\left(\frac{1}{n} \log \Vert A^{(n)}(\omega)v \Vert\right) = L_1(A,K)
	\end{equation}
for every $(\omega_0,\hat{v}) \in \Sigma\times\mathbb{S}^{d-1}$.

Assume first by contradiction that fixing $\omega_0\in \Sigma$, the convergence in $\hat{v}$ is not uniform, then there exist $\delta>0$ and, for simplicity of notation, a sequence $\{v_{n}\}_{n\in \N} \subset \mathbb{S}^{d-1}$ such that
$$\left|\mathbb{E}_{\omega_0}\left(\frac{1}{n} \log \Vert A^{(n)}(\omega)v_{n} \Vert\right) - L_1(A,K)\right| \geq \delta, \forall k\geq 1.$$

Since
\begin{align*}
\mathbb{E}_{\omega_0} \left(\frac{1}{n} \log \Vert A^{(n)}(\omega)v_{n} \Vert\right) &\leq \mathbb{E}_{\omega_0} \left(\frac{1}{n} \log \Vert A^{(n)}(\omega) \Vert\right)\\
& \to L_1(A,K) < L_1(A,K) + \frac{\delta}{2},
\end{align*}
we have that for $n\leq N$ with $N$ large enough, it can not happen that
$$\mathbb{E}_{\omega_0}\left(\frac{1}{n} \log \Vert A^{(n)}(\omega)v_{n} \Vert\right) \geq L_1(A,K) + \delta.$$
Thus we only need to consider the case with
$$\mathbb{E}_{\omega_0}\left(\frac{1}{n} \log \Vert A^{(n)}(\omega)v_{n} \Vert\right) \leq L_1(A,K) - \delta$$
We are going to prove that it cannot happen either. First, we claim that
$$\liminf_{n\to \infty} \frac{\Vert A^n(\omega)v_n \Vert}{\Vert A^n(\omega) \Vert} = c(\omega) > 0$$
for $\Pp_{\mu}$-almost every $\omega \in \Sigma^{\N}$. Note that

$$\frac{\Vert A^n(\omega)v_n \Vert}{\Vert A^n(\omega) \Vert} \geq |v_n \cdot \bar{v}^{(n)}(A)| \to |v \cdot \bar{v}^{(\infty)}(A)|,$$
where $\bar{v}^{(n)}(A)$ is the most expanding direction of each $n$-th iterates and $\bar{v}^{(\infty)}(A)$ is such that $\bar{v}^{(n)}(A) \to \bar{v}^{(\infty)}(A)$.

On the other hand, $\bar{v}^{(\infty)}(A)^{\bot}$ is the sum of all invariant subspaces in Osedelets decomposition associated with Lyapunov exponents $< L_1(A,K)$. 
Then, the quasi-irreducibility implies that
$$\liminf\frac{\Vert A^n(\omega)v_n \Vert}{\Vert A^n(\omega) \Vert} > 0$$
for $\Pp_{\mu}$-almost surely. Therefore,
$$\frac{1}{n}\log\frac{\Vert A^n(\omega)v_n \Vert}{\Vert A^n(\omega) \Vert} \to 0$$
almost surely as $n\to \infty$. Using the Dominated Convergence Theorem
\begin{align*}
	&\lim \frac{1}{n}\mathbb{E}_{\omega_0} \left[\log\Vert A^n(\omega)v_n \Vert \right]\\
	&= \lim \frac{1}{n}\mathbb{E}_{\omega_0} \left[\log\Vert A^n(\omega) \Vert \right] + \lim \frac{1}{n}\mathbb{E}_{\omega_0} \left[\log \frac{\Vert A^n(\omega)v_n \Vert}{\Vert A^n(\omega) \Vert} \right]\\
	&= L_1(A,K) + 0 = L_1(A,K)
\end{align*}
which establishes the claim and proves that the limit
$$\lim_{n\to\infty} \frac{1}{n} \mathbb{E}_{\omega_0} (\log \Vert A^{(n)}(\omega)v \Vert) = L_1(A,K),$$
is uniform convergence in $\hat{v}\in \mathbb{S}^{d-1}$ when $\omega_0\in \Sigma$ is fixed.

Now, it remains to prove that the previous limit converges uniformly in both $(\omega_0,\hat{v})$.

For every $M \in \GL_d(\R)$, consider
	$$
	\chi(M) = \sup\{\log \Vert M \Vert, \log \Vert M \Vert^{-1}\}.
	$$
and define $g_n: \Sigma^{\N} \to \R$ such that
	\begin{align*}
		& g_n(\omega) = g_n(\omega_0) \\
		=& \sup \left\{\left|\frac{1}{n} \mathbb{E}_{\omega_0} (\log \Vert A^{(n)}(\omega)v \Vert) - L_1(A,K)\right|: v\in \mathbb{S}^{d-1}, \omega = \{\omega_n\}\in \Sigma^{\N} \right\}	
	\end{align*}

Since $-\chi(M) \leq \log \Vert Mu \Vert \leq \chi(M)$ for every $u \in \mathbb{S}^{d-1}$, the sequence $g_n$ is uniformly bounded and
	$$
	|g_n(\omega)| \leq \frac{1}{n} \mathbb{E}_{\omega_0} (\chi(A^{(n)}(\omega)) + |L_1(A,K)|
	$$

Moreover, by (\ref{pointconver}), $g_n(\omega) \to 0$ when $n \to \infty$ for $\Pp_{\mu}$-a.e. $\omega \in \Sigma^{\N}$. Hence,
	$$
	\lim_{n\to \infty} \int g_n(\omega) \; d\Pp_{\mu}(\omega) = 0
	$$
Let $p \in \N$ such that $p < n$ and consider $v_p = \frac{A^p(\omega)v}{\Vert A^p(\omega)v \Vert}$,
	\begin{align*}
		&\Big|\frac{1}{n} \mathbb{E}_{\omega_0}  (\log \Vert A^{(n)}(\omega)v \Vert) - L_1(A,K)\Big| \leq \\
		& \leq \left|\frac{1}{n} \mathbb{E}_{\omega_0} (\log \Vert A^{(n-p)}(f^{p}\omega)v_p \Vert) - L_1(A,K)\right| + \left|\frac{1}{n} \mathbb{E}_{\omega_0} (\log \Vert A^{(p)}(\omega)v \Vert)\right| \\
		& \leq \left|\frac{1}{n} \mathbb{E}_{\omega_0} (\mathbb{E}_{\omega_p} (\log \Vert A^{(n-p)}(f^{p}\omega)v_p \Vert)) - L_1(A,K)\right| + \frac{1}{n} \mathbb{E}_{\omega_0} (\chi( A^{(p)}(\omega))) \\
		& \leq \mathbb{E}_{\omega_0}\left(\frac{1}{n-p} \mathbb{E}_{\omega_p} (\log \Vert A^{(n-p)}(f^{p}\omega)v_p \Vert) - L_1(A,K)\right) + \frac{p}{n} (|L_1(A,K)| + a)
	\end{align*}
Hence,
	$$
	g_n(\omega) \leq \mathbb{E}_{\omega_0}(g_{n-p}(f^p\omega)) + \frac{p}{n}(|L_1(A,K)| + a)
	$$

Since $K$ is uniformly ergodic, there exists a sequence $\epsilon(p)$, where $\epsilon(p) \to 0$ when $p\to \infty$, such that
	
	$$
	\sup_{\omega_0 \in \Sigma} \left|\int g_{n-p}(y) \; dK^p_{\omega_0}(y) - \int g_{n-p}(y) \; d\Pp_{\omega_0}(y)\right| \leq \epsilon(p)
	$$
for every $n$. Then, for every $p$,
	$$
	\lim_{n\to \infty} \sup_{\omega_0 \in \Sigma} g_n(\omega_0) \leq \lim_{n\to \infty} \int g_{n-p}(y) \; dK^p_{\omega_0}(y) + \epsilon(p) \leq \epsilon(p),
	$$
and this concludes the proof.
\end{proof}

\section{The strong mixing of the Markov operator}\label{unique}
Assuming the quasi-irreducibility of the cocycle $(A, K)$ and the simplicity of its top Lyapunov exponent, that is,  $L_1(A,K)>L_2(A,K)$, we prove that the powers of the associated Markov operator $\bar Q^n$ converge in an appropriate sense to its (eventually unique) stationary measure. We also derive a version of Furstenberg's formula that will then allow us to prove the H\"older continuity of the Lyapunov exponents. 
We follow closely the approach  in~\cite[Section 5.3.2]{DK-book}.

\subsection{Contracting property of the H\"older seminorm}
We show that the Markov operator acts as a contraction on an appropriate space of observables.
 
Consider on the projective space $\Proj$ the distance
$$
\delta(\hat{p},\hat{q}):= \frac{\| p \wedge q \|}{\|p\| \|q\|},
$$
where $p$ and $q$ are representatives of $\hat{p}$ and $\hat{q}$ respectively.

Given $0 < \alpha \leq 1 $ and $\varphi \in L^\infty(\Sigma \times \Proj)$, we define the H\"older seminorm $v_\alpha$, the H\"older norm $\| \cdot \|_\alpha$ and the space of H\"older continuous functions ${\mathcal{H}}_\alpha(\Sigma \times \Proj)$ by:
\begin{align*}
&v_\alpha = \sup_{ \substack  {{\omega_0 \in \Sigma} \\    \hat{p} \neq \hat{q}}  } \frac{| \varphi(\omega_0,\hat{p}) - \varphi(\omega_0,\hat{q}) |}{\delta(\hat{p},\hat{q})^\alpha}, \\
&\| \varphi \|_\alpha = v_\alpha(\varphi) + \| \varphi \|_\infty, \\
&{\mathcal{H}}_\alpha(\Sigma \times \Proj) = \{ \varphi \in L^\infty(\Sigma \times \Proj) \colon \|\varphi\|_\alpha < \infty \}.
\end{align*}

Moreover, consider the average H\"older constant:
$$
k_\alpha(A,K) = \sup_{\substack  {{\omega_0 \in \Sigma} \\    \hat{p} \neq \hat{q}}} \int \frac{\delta(\hat{A}(\omega_1,\omega_0)\hat{p}, \hat{A}(\omega_1,\omega_0)\hat{q})^\alpha}{\delta(\hat{p},\hat{q})^\alpha} \;dK_{\omega_0}(\omega_1).
$$

\begin{proposition} \label{contracaoValpha}
For all $\varphi \in {\mathcal{H}}_\alpha(\Sigma \times \Proj)$, 
$$
v_\alpha (Q_{A,K} (\varphi)) \leq v_\alpha(\varphi)k_\alpha(A,K).
$$
\end{proposition}

\begin{proof}
\begin{align*}
v_\alpha(Q_{A,K}(\varphi)) &= \sup_{\substack  {{\omega_0 \in \Sigma} \\    \hat{p} \neq \hat{q}}} \frac{\left | \int \varphi (\omega_1, \hat{A}(\omega_1,\omega_0)\hat{p}) -  \varphi(\omega_1,\hat{A}(\omega_1,\omega_0)\hat{q})\;dK_{\omega_0}(\omega_1) \right|}{\delta(\hat{p},\hat{q})^\alpha}   \\
&\leq \sup_{\substack  {{\omega_0 \in \Sigma} \\    \hat{p} \neq \hat{q}}} \frac{\int \left |  \varphi (\omega_1, \hat{A}(\omega_1,\omega_0)\hat{p}) -  \varphi(\omega_1,\hat{A}(\omega_1,\omega_0)\hat{q}) \right| \;dK_{\omega_0}(\omega_1)}{\delta(\hat{p},\hat{q})^\alpha}  \\
&\leq v_\alpha(\varphi) \sup_{\substack  {{\omega_0 \in \Sigma} \\    \hat{p} \neq \hat{q}}} \int \frac{\delta(\hat{A}(\omega_1,\omega_0)\hat{p}, \hat{A}(\omega_1,\omega_0)\hat{q})^\alpha}{\delta(\hat{p},\hat{q})^\alpha} \;dK_{\omega_0}(\omega_1) \\
&\leq v_\alpha(\varphi) k_\alpha(A,K).
\end{align*}
\end{proof}

\begin{proposition} \label{k_alpha submultiplicative}
The sequence $\{k_\alpha(A^n,K^n)\}_n$ is sub-multiplicative:
$$
k_\alpha(A^{m+n},K^{m+n}) \leq k_\alpha(A^n,K^n)k_\alpha(A^m,K^m).
$$
\end{proposition}

\begin{proof}By definition,
\begin{align*}
k_\alpha(A^{m+n},K^{m+n})= \sup_{\substack  {{\omega_0 \in \Sigma} \\    \hat{p} \neq \hat{q}}} \int \frac{\delta(\hat{A}^{m+n}(\omega)\hat{p}, \hat{A}^{m+n}(\omega)\hat{q})^\alpha}{\delta(\hat{p},\hat{q})^\alpha} \;dK^{m+n}_{\omega_0}(\omega_1,\omega_2, \dots , \omega_{m+n}).
\end{align*}
Multiply and divide inside the integral by $\delta(\hat{A}^{m}(\omega_1,\omega_0)\hat{p}, \hat{A}^{m}(\omega_1,\omega_0)\hat{q})^\alpha$ and conclude that it is less or equal to $k_\alpha(A^n,K^n)k_\alpha(A^m,K^m)$.
\end{proof}

Recall that the space of cocycles
\begin{align*}
\mathcal{C} := \{ (A,K)\colon &A \colon \Sigma \times \Sigma \to \GL_d(\R) \text{ is Lipschitz continuous and } \\
 &K \colon \Sigma \to \Prob(\Sigma) \text{ is uniformly ergodic and}\\ &\text{continuous in the weak* topology} \}
\end{align*}
is endowed with the metric
$$
d( (A,K), (B,L) ) = \max \{ d_\infty(A,B), d_{W_1}(K, L) \} ,
$$
where the distance between two Markov kernels is given by:
$$
d_{W_1}(K, L) = \sup_{\om_0 \in \Sigma} W_1 (K_{\om_0}, L_{\om_0}) \,.
$$
$W_1$ is the Wasserstein distance, which metrizes the weak* topology on $\Prob (\Sigma)$.  
The Kantorovich-Rubinstein theorem characterizes the Wasserstein distance as follows:
\begin{align}
W_1(\mu, \nu)&= \inf_{\pi \in \Pi(\mu,\nu)} \int d(x,y) \; d \pi \\
&= \sup_{\varphi \; \in \Lip_1 (\Sigma)} \ \  \int \varphi \; d\mu - \int \varphi \; d\nu ,
\end{align}
where $\Pi(\mu,\nu)$ is the set of all coupling measures of $\mu$ and $\nu$ (probability measures on $\Sigma \times \Sigma$ with marginals $\mu$ and $\nu$) and $\Lip_1 (\Sigma)$ is the set of Lipschitz continuous functions on $\Sigma$ with Lipschitz constant $\le 1$.

\begin{proposition} \label{Lipp}
Fix $n \in \N$. The map $(A,K) \mapsto (A^n,K^n)$ is Lipschitz with respect to the metric $d$.
\end{proposition}

\begin{proof}
Note that the linear map $A \mapsto A^n$ is Lipschitz with constant $C(n)$ that depends on $n$ but not on the kernel. Moreover, we claim that the map $K\mapsto K^n$ is also Lipschitz with constant $n$:
\begin{align*}
&d(K^2,L^2) = \sup_{\substack  {{\omega_0 \in \Sigma} \\    \varphi \in \Lip_1(\Sigma)}} \left | \int \varphi(\omega_2) \; d(K^2_{\omega_0}(\omega_2) - L^2_{\omega_0}(\omega_2))  \right| \\
&\leq \sup_{\substack  {{\omega_0 \in \Sigma} \\    \varphi \in \Lip_1(\Sigma)}}   \left| \int \varphi(\omega_2) \; dK_{\omega_1}(\omega_2)dK_{\omega_0}(\omega_1 ) -  \int \varphi(\omega_2) \; dK_{\omega_1}(\omega_2)dL_{\omega_0}(\omega_1) \right| + \\
&+ \left| \int \varphi(\omega_2) \; dK_{\omega_1}(\omega_2)dL_{\omega_0}(\omega_1 ) -  \int \varphi(\omega_2) \; dL_{\omega_1}(\omega_2)dL_{\omega_0}(\omega_1) \right| \\
&\leq 2d(K,L).
\end{align*} 
The proof of the claim follows by induction. Hence, the joint map $(A,K) \mapsto (A^n,K^n)$ is Lipschitz with constant the maximum between $C(n)$ and $n$.
\end{proof}

\begin{proposition} \label{Kalphamenorque1}
Let $(A,K) \in (\mathcal{C},d)$. Assume that 
\begin{enumerate}[(i)]
\item $A$ is quasi irreducible with respect to $(K,\mu)$,
\item $L_1(A,K) > L_2(A,K)$.
\end{enumerate}
Then, there are numbers $\delta >0$, $0<\alpha<1$, $0< \sigma <1$ and $n \in \N$ such that for all $(B,L) \in (\mathcal{C},d)$ with $d((B,L),(A,K)) < \delta$ one has $k_\alpha(B^n,L^n) < \sigma $.
\end{proposition}

\begin{proof}
We begin the proof with the following lemma:

%
%

\begin{lemma} 
Given a pair $(A,K) \in (\mathcal{C},d)$, then for all $\alpha>0$,
$$
k_\alpha(A,K) \leq \sup_{\substack  {{\omega_0 \in \Sigma} \\    \hat{p} \in \Proj}} \int_\Sigma \left({\frac{|s_1( \hat{A}(\omega_1,\omega_0)) s_2( \hat{A}(\omega_1,\omega_0))|}{\|\hat{A}(\omega_1,\omega_0) \hat{p}\|^2}}\right)^\alpha \;dK_{\omega_0}(\omega_1) ,
$$
where 
$s_1(\cdot)$ and $s_2(\cdot)$ are the first and second singular values.
\end{lemma}

\begin{proof} 

We claim that given $\alpha > 0$ and two points $\hat{p}, \hat{q} \in \Proj$, we have 
$$
\left[ \frac{\delta(\hat{A}(\omega) \hat{p}, \hat{A}(\omega) \hat{q})}{\delta(\hat{p}, \hat{q})} \right] ^\alpha \leq \frac{| s_1(A(\omega))s_2(A(\omega))|^\alpha}{2} \left[ \frac{1}{\| A(\omega)p  \|^{2\alpha}} + \frac{1}{\| A(\omega)q  \|^{2\alpha}} \right]
$$
for every $\omega_0 \in \Sigma$. 

Note that if we integrate with respect to the measure $K_{\omega_0}$ and take the supremum in $\hat{p} \neq \hat{q}$ on both sides of this inequality, we conclude the lemma. Therefore it is enough to prove the previous claim.

By the exterior product property, 
$$
\| A(\omega_1, \omega_0)p \wedge A(\omega_1, \omega_0)q \| = |s_1(A(\omega_1, \omega_0))s_2(A(\omega_1, \omega_0))|\|p \wedge q\|.
$$
Hence, by the definition of the projective distance and the fact that the geometric mean is less or equal the arithmetic mean,
\begin{align*}
\left[ \frac{\delta(\hat{A}(\omega) \hat{p}, \hat{A}(\omega) \hat{q})}{\delta(\hat{p}, \hat{q})} \right] ^\alpha &= \left[ \frac{|s_1(A(\omega))s_2(A(\omega))|}{\|A(\omega)p\| \|A(\omega)q\|} \right]^\alpha \\
&\leq \frac{|s_1(A(\omega))s_2(A(\omega))|^\alpha}{2} \left[ \frac{1}{\|A(\omega)p\|^{2\alpha}} +\frac{1}{\|A(\omega)q\|^{2\alpha}}  \right].
\end{align*}

\end{proof}

Now we proceed with the proof of the proposition. 

By theorem $\ref{convUnif}$, given $(A,K) \in (\mathcal{C},d)$ 
satisfying assumptions $i$ and $ii$, we have
$$\lim_{n\to\infty} \frac{1}{n} \mathbb{E}_{\omega_0} (\log \Vert A^{(n)}(\omega)v \Vert ^{-2}) =-2 L_1(A,K),$$
with uniform convergence in $(w_0,\hat{v})\in \Sigma\times \mathbb{S}^{d-1}$. 

Hence, by choosing $\epsilon$ small enough e.g $\frac{1}{4}(L_1(A,K) - L_2(A,K))$ and $n$ sufficiently large, we conclude that
$$
\mathbb{E}_{\omega_0} (\log \Vert A^{(n)}(\omega)v \Vert ^{-2}) \leq n (-2L_1(A,K) +\epsilon) .
$$
Moreover, we know that
$$
\log|s_1(A^{(n)}(\omega))|+\log|s_2(A^{(n)}(\omega))|\leq n(L_1(A,K) + L_2(A,K) + \epsilon).
$$

We combine these two estimates to conclude that, for $n$ sufficiently large, 
\begin{equation}\label{compute eq}
 \mathbb{E}_{\omega_0} \log \left[ \frac{s_1( A^{(n)}(\omega)) s_2( A^{(n)}(\omega))}{ \Vert A^{(n)}(\omega)v \Vert ^{2}} \right] \leq -1 ,
\end{equation}
since $L_1(A,K)>L_2(A,K)$.

By the inequality 
$
e^x \leq 1+x+\frac{x^2}{2}e^{|x|},
$
we conclude that for every $v \in \mathbb{S}^{d-1}$ and every $\omega_0 \in \Sigma$, 
\begin{align*}
& \mathbb{E}_{\omega_0} \exp \left( \log \left[ \frac{| s_1( A^{(n)}(w)) s_2( A^{(n)}(\omega))|}{ \Vert A^{(n)}(\omega)v \Vert ^{2}} \right]^\alpha \right) \leq  \\
& \leq   1+\mathbb{E}_{\omega_0} \left[ \alpha \log \frac{|s_1( A^{(n)}(\omega)) s_2( A^{(n)}(\omega))|}{ \Vert A^{(n)}(\omega)v \Vert ^{2}} \right] + \\
&+\mathbb{E}_{\omega_0} \left[ \frac{\alpha^2}{2}\log^2  \frac{|s_1( A^{(n)}(\omega)) s_2( A^{(n)}(\omega))|}{ \Vert A^{(n)}(\omega)v \Vert ^{2}}e^{\frac{| \alpha \log | s_1( A^{(n)}(\omega)) s_2( A^{(n)}(\omega))||}{ \Vert A^{(n)}(\omega)v \Vert ^{2}} } \right] \\
&\leq 1- \alpha +C \frac{\alpha^2}{2}.
\end{align*}
Hence, $k_\alpha(A^n,K^n) \leq 1- \alpha +C \frac{\alpha^2}{2}$.
Note that $C$ is a constant that depends only on $(A,K)$ and $n$. Thus, we can choose $\alpha$ small enough such that $k_\alpha(A^n,K^n)<1$. Moreover, $k_\alpha(A,K)$ depends continuously on $(A,K)$ and, by proposition \ref{Lipp}, the map $(A,K) \mapsto  (A^{(n)},K^n)$ is Lipschitz, therefore we can extend the result to a neighborhood of $(A,K)$.
\end{proof}

\subsection{The strong mixing property} We begin with the general concept of strong mixing. 
Let $(M, K, \mu)$ be a Markov system. We need to define properly the space of observables on which the Markov operator $Q$ associated with $K$ acts. Let $(\Escr, \norm{\cdot}_\Escr)$ be a Banach space where $\Escr \subset C^0(M)$ is $Q$-invariant in the sense that $\varphi\in\Escr \Leftrightarrow Q\varphi \in \Escr$. Moreover, we assume that the constant function $\ind \in \Escr$ and that the inclusion of $\Escr \subset C^0(M)$ is continuous, namely $\norm{\varphi}_\infty \leq C_1 \norm{\varphi}_\Escr$ for some constant $C_1<\infty$. We also assume that $Q$ is bounded (or continuous) on $(\Escr, \norm{\cdot}_\Escr)$, i.e. $\norm{Q\varphi}_\Escr \leq C_2 \norm{\varphi}_\Escr$ with $C_2 <\infty$. In practice we will have $C_1=C_2=1$.

\begin{definition}[Strong mixing]
	The Markov system $(M,K,\mu,\Escr)$ (or simply $Q$) is strongly mixing if there are $C<\infty$ and $\sigma \in (0, 1)$ such that for all 
	$n\in \N$ and $\varphi\in \Escr$,
	$$
	\norm{Q^n\varphi-\int_M \varphi d\mu}_\infty \leq C \sigma^n \norm{\varphi}_\Escr \, .
	$$
\end{definition}

By assumption, we have $K_{\omega_0}^n \to \mu$ uniformly in $\omega_0\in \Sigma$. By Remark \ref{equivalence}, $\forall \,\varphi \in L^\infty(\mu)$, 
$$
\norm{Q^n\varphi-\int \varphi d\mu}_\infty \leq C\sigma^n \norm{\varphi}_\infty
$$
for some constant $C<\infty$ and $0<\sigma<1$.

With the contracting property of the $v_\alpha$ seminorm in hand, we are ready to prove that our Markov operator $Q_A$ is strongly mixing on the space of H\"older functions $\mathcal{H}_\alpha(\Sigma\times \Proj)$. Indeed, this holds in a neighborhood of $A$ as shown by the following theorem.

\begin{theorem}\label{teoStrongMixing} 
	Given $(A,K) \in (\mathcal{C},d)$ such that the assumptions of proposition \ref{Kalphamenorque1} are satisfied, then there exist constants $C<\infty$, $0<\sigma<1$ and a neighborhood $U$ of $(A,K)$ in $ (\mathcal{C},d)$ such that for all $(B,L) \in U$, $Q_{B,L}$ is strongly mixing on $\mathcal{H}_\alpha(\Sigma\times \Proj)$:
	$$
	\norm{Q_{B,L}^{n}\varphi-\int_{\Sigma\times \Proj} \varphi d\eta_{B,L}}_\infty \leq C \sigma^n \norm{\varphi}_\alpha, \, \forall \, \varphi\in \mathcal{H}_\alpha(\Sigma\times \Proj).
	$$
	Moreover, since the $v_\alpha$ seminorm is also exponentially contracting on $\mathcal{H}_\alpha(\Sigma\times \Proj)$, we further get
	$$
	\norm{Q_{B,L}^{n}\varphi-\int_{\Sigma\times \Proj} \varphi d\eta_{B,L}}_\alpha \leq C \sigma^n \norm{\varphi}_\alpha, \, \forall \, \varphi\in \mathcal{H}_\alpha(\Sigma\times \Proj).
	$$
	In fact, $\sigma$ can be choosen as the square root of the minimum between the contracting rate of the $v_\alpha$ seminorm and the convergence rate of the kernel $L: \Sigma \to \Prob(\Sigma)$, and we may always choose them to be the same by setting properly the size of the neighbourhood.
	\end{theorem}

\begin{proof}
	Take $(B,L)\in U$ where $U$ is given by Proposition \ref{Kalphamenorque1}, $n\geq m$ with $n,m\in \N$, $\varphi\in \mathcal{H}_\alpha(\Sigma\times \Proj)$ and $\eta_{B,L}$ any $L_B$-stationary measure, define 4 families of transformations in the following way: for any $(\omega_0, p)\in \Sigma \times \Proj$,
	\begin{enumerate}
		\item $(T_{B,L,n}^{(0)}\varphi)(\omega_0,p):=(Q_{B,L}^n \varphi)(\omega_0,p)=\EE_{\omega_0}[\varphi(\omega_n, B^{(n)}p)]$.
		\item $(T_{B,L,n,m}^{(1)}\varphi)(\omega_0,p):=\EE_{\omega_0}[\varphi(\omega_n, (B^{(m)}\circ T^{n-m})p)]$.
		\item $(T_{B,L,m}^{(2)}\varphi)(\omega_0,p):=\EE_\mu[\varphi(\omega_n,B^{(m)}p)]$, constant in $\omega_0$, thus we denote it by $(T_{B,L,m}^{(2)}\varphi)(p)$ which is a compact transformation.
		\item $(T_{B,L}^{(3)}\varphi)(\omega_0,p):=\int \varphi d\eta_{B,L}$, constant.
	\end{enumerate}
Then it is straightforward to obtain the following inequalities:
\begin{enumerate}
	\item $\abs{(T_{B,L,n}^{(0)}\varphi)(\omega_0,p)-T_{B,L,n,m}^{(1)}\varphi)(\omega_0,p)} \leq C \sigma^m \norm{\varphi}_\alpha$ for the same $\sigma$ using the contracting property of the $v_\alpha$ seminorm.
	\item $\abs{T_{B,L,n,m}^{(1)}\varphi)(\omega_0,p)-(T_{B,L,m}^{(2)}\varphi)(p)}\leq C \sigma^{n-m}\norm{\varphi}_\alpha$ using the uniform convergence rate of $L_{\omega_0}^{n-m} \to \mu$.
	\item $\abs{(T_{B,L,m}^{(2)}\varphi)(p)-(T_{B,L,n}^{(2)}\varphi)(p)}\leq C\sigma^m \norm{\varphi}_\alpha$ using again the contracting property of the $v_\alpha$ seminorm.
\end{enumerate}

For simplicity, we may set $n=2m$ in $(1)$ and $(2)$, and set $n=l$ in $(3)$, then by $(1)$-$(3)$, we have for all $B\in U$ and $\varphi \in\mathcal{H}_\alpha(\Sigma\times \Proj)$,
$$
\norm{Q_{B,L}^{2m}\varphi-T_{B,L,l}^{(2)}\varphi}_\infty \leq 3C\sigma^m \norm{\varphi}_\alpha.
$$
Note that the sequence $\{T_{B,L,l}^{(2)}\varphi\}_{l\geq0}$ is relatively compact in $C(\Proj)$. Then the set $S_\varphi$ of its limit points in $(C(\Proj), \norm{\cdot}_\infty)$ is non-empty. Take any $g\in S_\varphi$, we claim that 
$$
g=\int \varphi d\eta_{B,L}= T_{B,L}^{(3)}\varphi.
$$

Let us prove the claim. Take a subsequence $l_j \to \infty$ such that $\{T_{B,L,l_j}^{(2)}\varphi\}_{j\geq0}$ converges to $g$ in the previous inequality, we get
$$
\norm{Q_{B,L}^{2m}\varphi-g}_\infty \leq 3C\sigma^m\norm{\varphi}_\alpha.
$$

On the other hand, we have $v_\alpha(Q_{B,L}^{2m}\varphi)\leq C\sigma^{2m}\norm{\varphi}_\alpha$. This implies $v_\alpha(g)=0$ which further implies that $g$ is constant. Finally, using the condition that $\eta_{B,L}$ is $L_B$-stationary, we have
$$
\int Q_{B,L}^{2m}\varphi d\eta_{B,L}=\int \varphi d\eta_{B,L}
$$
which equals $g$. Take $\sigma'=\sigma^{\frac{1}{2}}$ as a new parameter, this finishes the whole proof.
\end{proof}

\begin{corollary}[Uniqueness of the stationary measure] \label{uniquenessofStatMeas} 
	Given $(A,K) \in (\mathcal{C},d)$ such that the assumptions of proposition \ref{Kalphamenorque1} are satisfied,  the kernel $L_B$ on $\Sigma \times \Proj$ has a unique stationary measure $\eta_{B,L}$ for every $(B,L)\in U$, which further gives that $\bar{L}_B$ has a unique stationary measure $L \ltimes \eta_{B,L}$.
\end{corollary}

\begin{proof}
	Assume there are two different stationary measures $\eta_{B,L}$ and $\eta'_{B,L}$, using Theorem \ref{teoStrongMixing} and the fact that $\mathcal{H}_\alpha(\Sigma \times \Proj)$ is dense in $L^\infty(\Sigma \times \Proj)$, we get $\eta_{B,L}=\eta'_{B,L}$.
\end{proof}

As a result, we can finally upgrade the Furstenberg formula as follows.

\begin{theorem} \label{FurstFormula}
	Given $(A,K) \in (\mathcal{C},d)$ such that the assumptions of proposition \ref{Kalphamenorque1} are satisfied, there exists a neighborhood $U$ of $(A,K)\in(\mathcal{C},d)$ such that for every $(B,L)\in U$
	$$
	L_1(B,L)= \int_{\Sigma \times \Sigma \times \Proj} \psi(y,x,\hat{v})\; dL_x(y)d\eta_{B,L}(x,\hat{v})
	$$
	where $\psi: \Sigma \times \Sigma \times \Proj \to \R$ is such that
	\begin{equation} \label{observableLyapunov}
		\psi(y,x,\hat{v}) = \log \frac{\Vert A(y,x)v \Vert}{\Vert v \Vert}
	\end{equation}
and $\eta_{B,L}$ is the unique $L_B$-stationary measure.
\end{theorem}



\section{Continuity of the Lyapunov exponent} \label{joint}
In this section we prove the joint $(A,K)\mapsto L_1(A,K)$ H\"older continuity of the Lyapunov exponents of random Markov cocycles based on the technique introduced by Baraviera and Duarte in \cite{BD}. This approach, via the Furstenberg Formula, enables the study of this type of continuity without the need of going through the theory of large deviations. Moreover, it also has the advantage of providing a computable H\"older exponent.

We start with an adaptation of Proposition 4.7 of \cite{DK-31CBM}, which shows that for every random cocycle satisfying generic hypothesis and every $n \in \N$, the map $A \mapsto Q^n_A$ is locally H\"older. We extend it to the mixing Markov case in $(A,K)$. More precisely,

\begin{proposition}
Let $(A,K) \in (\mathcal{C},d).$ 
Assume that:

\begin{enumerate}[(i)]
\item $A$ is quasi irreducible with respect to $(K,\mu)$,
\item $L_1(A,K) > L_2(A,K)$.
\end{enumerate}

Then, there exists $\delta > 0 $, such that for all $(B,L)$ and $(D,T)$ in $(\mathcal{C},d)$
satisfying $d((A,K),(B,L)) < \delta$ and $d((A,K),(D,T)) < \delta$, for all $f \in \mathcal{H}_\alpha(\Sigma \times \Proj)$ and every $n \in \N$,
$$
\|Q_{B,L}^nf - Q_{D,T} ^n f \|_\infty \leq C d((B,L),(D,T))^\alpha v_\alpha(f).
$$
\end{proposition}

\begin{proof}
First consider the case $n=1$. For every $\omega_0 \in \Sigma$ and $v \in \Proj$,
\begin{align*}
&\left\| (Q_{B,L} - Q_{D,T}) (f) \right\|_\infty  = \\
&= \sup_{ \substack{ \hat{v} \in \Proj \\ \omega_0 \in \Sigma} } \left| \int_\Sigma f (\omega_1, B(\omega_1,\omega_0)v)\; dL_{\omega_0}(\omega_1) - \int_\Sigma f(\omega_1,D(\omega_1,\omega_0)v) \; dT_{\omega_0}(\omega_1) \right| \\
&\leq \sup_{ \substack{ \hat{v} \in \Proj \\ \omega_0 \in \Sigma} } \left| \int_\Sigma f (\omega_1, B(\omega_1,\omega_0)v)- f(\omega_1,D(\omega_1,\omega_0)v) \; dL_{\omega_0}(\omega_1) \right| + \\
& + \sup_{ \substack{ \hat{v} \in \Proj \\ \omega_0 \in \Sigma} } \left| \int_\Sigma f (\omega_1, D(\omega_1,\omega_0)v)\; dL_{\omega_0}(\omega_1) - \int_\Sigma f(\omega_1,D(\omega_1,\omega_0)v) \; dT_{\omega_0}(\omega_1) \right|.
\end{align*}
Since $f$ is H\"older, we can bound the first term by
\begin{align*}
&\sup_{ \substack{ \hat{v} \in \Proj \\ \omega_0 \in \Sigma} } \left| \int_\Sigma f (\omega_1, B(\omega_1,\omega_0)v)- f(\omega_1,D(\omega_1,\omega_0)v) \; dL_{\omega_0}(\omega_1) \right| \leq \\
&\leq \int_\Sigma v_\alpha(f) \; \delta(B(\omega_1,\omega_0)v,D(\omega_1,\omega_0)v)^\alpha \; dL_{\omega_0}(\omega_1) \\
&\leq v_\alpha(f)\; d_\infty(B,D)^\alpha .
\end{align*}

Now we proceed to estimate the second term. For every $\pi \in \Pi(L_{\omega_0},T_{\omega_0})$:
\begin{align*}
& \sup_{ \substack{ \hat{v} \in \Proj \\ \omega_0 \in \Sigma} }  \left| \int_\Sigma f (\omega_1, D(\omega_1,\omega_0)v) \; dL_{\omega_0}(\omega_1)- \int_\Sigma f (z_1, D(z_1,\omega_0)v) \; dT_{\omega_0}(z_1) \right| \\
&= \sup_{ \substack{ \hat{v} \in \Proj \\ \omega_0 \in \Sigma} } \left|  \int_{\Sigma \times \Sigma}  f (\omega_1,D(\omega_1,\omega_0)v)-f (z_1,D(z_1,\omega_0)v) \; d\pi(\omega_1,z_1) \right|  \\
&\leq v_\alpha(f) \sup_{ \substack{ \hat{v} \in \Proj \\ \omega_0 \in \Sigma} } \int_{\Sigma \times \Sigma} \delta \left(D(\omega_1,\omega_0)v , D(z_1,\omega_0)v \right)^\alpha \; d\pi(\omega_1,z_1)  \\
&\leq v_\alpha(f) \sup_{ \substack{ \hat{v} \in \Proj \\ \omega_0 \in \Sigma} } \int_{\Sigma \times \Sigma} \|D(\omega_1,\omega_0) - D(z_1,\omega_0)\|^\alpha \times \\
&\times \max \left \{ \frac{1}{\|D(\omega_1,\omega_0)(v)\|}, \frac{1}{\|D(z_1,\omega_0)(v)\|} \right\}^\alpha \; d\pi(\omega_1,z_1).
\end{align*}

Since $\Sigma$ is compact, there exists a constant $C_1>0$, such that $\max \left \{ \frac{1}{\|D(\omega_1,\omega_0)(v)\|}, \frac{1}{\|D(z_1,\omega_0)(v)\|} \right\}^\alpha \leq  C_1$. 
Then, for every $  \pi \in \Pi(L_{\omega_0},T_{\omega_0}) $, we can bound the second term by
\begin{align*}
& C_1  v_\alpha(f) \sup_{\omega_0 \in \Sigma}   \int_{\Sigma \times \Sigma} \|D(\omega_1,\omega_0) - D(z_1,\omega_0)\|^\alpha \; d\pi(\omega_1,z_1) \\
&\leq  C_2  v_\alpha(f) \sup_{\omega_0 \in \Sigma } \left(  \int_{\Sigma \times \Sigma} \|(\omega_1,\omega_0) - (z_1,\omega_0)\|  \; d\pi(\omega_1,z_1) \right)^\alpha \\
&\leq  C_2  v_\alpha(f) \; d(L, T)^\alpha,
\end{align*}
where on the second line we used the Lipschitz continuity of the map $x \mapsto D(x,y) $ and Jensen's inequality together with the concavity of the function $t \mapsto t^\alpha$, which holds when $t \in [0,\infty)$ and $\alpha \in (0,1]$.

Therefore, we conclude the case $n=1$: 
$$
\left\| (Q_{B,L} - Q_{D,T}) (f) \right\|_\infty \leq   C_2  v_\alpha(f) \; d((B,L), (D,T))^\alpha.
$$

Now we complete the proof using the following relation. Note that 
$$
Q_{B,L}^n  - Q_{D,T}^n = \sum_{i=0}^{n-1}Q^i_{D,T} \circ (Q_{B,L} - Q_{D,T}) \circ Q_{B,L}^{n-i-1}.
$$

Then, using the triangle inequality, the fact that the norm of the Markov operator is $1$ and the case $n=1$, we obtain:
\begin{align*}
\|Q_{B,L}^n(f)  - Q_{D,T}^n(f)\| &= \left \| \sum_{i=0}^{n-1} Q^i_{D,T} \circ (Q_{B,L} - Q_{D,T}) \circ Q_{B,L}^{n-i-1} (f) \right\|_\infty \\
&\leq \sum_{i=0}^{n-1} \left \| Q^i_{D,T} \circ (Q_{B,L} - Q_{D,T}) \circ Q_{B,L}^{n-i-1} (f) \right \|_\infty \\
&\leq \sum_{i=0}^{n-1} \left \| (Q_{B,L} - Q_{D,T}) \circ Q_{B,L}^{n-i-1} (f) \right \|_\infty \\
&\leq C_2  d((B,L),(D,T))^\alpha \sum_{i=0}^{n-1} v_\alpha\left( Q_{B,L}^{n-i-1}(f)\right).
\end{align*}
Since the operator contracts its seminorm $v_\alpha$ (see propositions \ref{contracaoValpha} and \ref{Kalphamenorque1}), we conclude that there exists $\delta>0$ and $\sigma < 1$ such that, if $d((A,K),(B,L))< \delta$, then $k_\alpha(B^n,L^n) < \sigma$. Moreover, since $k_\alpha$ is sub multiplicative, we conclude:
\begin{align*}
\|Q_{B,L}^n(f)  - Q_{D,T}^n(f)\| &\leq  C_2 d((B,L),(D,T))^\alpha v_\alpha (f) \sum_{i=0}^{\infty} k_\alpha(B^i,L^i)  \\
&\leq C d((B,L),(D,T))^\alpha v_\alpha(f).
\end{align*}

\end{proof}

Once the map $(A,K) \mapsto Q^n_{A,K}$ is locally H\"older and $Q^n_{A,K}$ converges to the stationary measure $\eta_{A,K}$ (in the sense of \ref{teoStrongMixing}), we now prove that the map $(A,K) \mapsto \eta_{A,K}$ is also locally H\"older.

\begin{corollary} \label{estimativa contin med est}
Given $(A,K) \in (\mathcal{C},d)$ such that the assumptions of proposition \ref{Kalphamenorque1} are satisfied, there exists $\delta >0$ such that for all $(B,L)$ and $(D,T)$ satisfying $d((A,K),(B,L))< \delta$ and $d((A,K),(D,T))< \delta$ and  for all $f \in \mathcal{H}_\alpha(\Sigma \times~ \Proj)$,
$$
\left | \int f \; d\eta_{B,L} - \int f \; d\eta_{D,T} \right | \leq C d((B,L),(D,T)) ^\alpha v_\alpha(f).
$$
\end{corollary}

\begin{proof}
By lemma $\ref{uniquenessofStatMeas}$, there are unique stationary measures $\eta_{B,L}$ and $\eta_{D,T}$ associated with the Markov kernels $L_B$ and $T_D$ respectively. Moreover, by theorem $\ref{teoStrongMixing}$,
%
$$
\lim_{n \to \infty} Q_{B,L}^n (f) = \left( \int f \; d\eta_{B,L} \right) \; \textbf{1} \quad \text{and} \quad \lim_{n \to \infty} Q_{D,T}^n (f) = \left( \int f \; d\eta_{D,T} \right) \; \textbf{1},
$$


where $1(\omega_0,v)=1$ for every $(\omega_0,v)$ is the constant function. Therefore, we conclude that:
\begin{align*}
\left | \int f \; d\eta_{B,L} - \int f \; d\eta_{D,T} \right | &\leq \sup_{n \to \infty} \| Q^n_{B,L}(f) - Q^n_{D,T}(f) \|\\
& \leq C d((B,L),(D,T))^\alpha v_\alpha(f).
\end{align*}
\end{proof}

An immediate consequence of the previous corollary is that the map $(A,K) \mapsto~ m_{A,K}:=~K \ltimes ~\eta_{A,K}$ is also locally H\"older, where $m_{A,K}$ is the unique stationary measure on $\Prob(\Sigma \times \Sigma \times \Proj)$ associated to the kernel $\bar{K}_A$.

We are now ready to prove the local H\"older continuity of the Lyapunov exponents.

\begin{theorem} \label{teorema principal}
Let $(A,K) \in (\mathcal{C},d)$. Assume that:
\begin{enumerate}[(i)]
\item $A$ is quasi irreducible with respect to $(K,\mu)$,
\item $L_1(A,K) > L_2(A,K)$.
\end{enumerate}
Then, there exists a neighbourhood $V$ of $(A,K)$ in $(\mathcal{C},d)$ where the map $(A,K) \mapsto L_1(A,K)$ is H\"older continuous.
\end{theorem}

\begin{proof}
By hypothesis $(i)$ and $(ii)$, we are in the setting of theorem $\ref{FurstFormula}$, thus we can express the top Lyapunov exponent $L_1(A,K)$ as
$$
L_1(A,K) = \int \psi_A(\omega_1,\omega_0,\hat{v})\; dK_{\omega_0}(\omega_1)d\eta_{A,K}(\omega_0,\hat{v})  =\int \psi_{A} \; dm_{A,K}  \;,
$$
where
$\psi_A(\omega_1,\omega_0,\hat{v}) = \log \frac{\Vert A(\omega_1,\omega_0)v \Vert}{\Vert v \Vert}$
and $m_{A,K}$ is the unique stationary measure associated to the Markov kernel $\bar{K}_A$.

Moreover, there exists a neighbourhood $V$ of $(A,K)$ in $(\mathcal{C},d)$ such that  for every $(B,L)$ and $(D,T)$ in $V$, we can express their top Lyapunov exponent using Furstenberg's Formula.

Therefore, by corollary \ref{estimativa contin med est} and the fact that $A \mapsto \psi_{A}$ is locally Lipschitz, we estimate:
\begin{align*}
&\left| L_1(B,L) - L_1(D,T) \right | = \left | \int \psi_{B} \; dm_{B,L} -  \int \psi_{D} \; dm_{D,T} \right | \\
&\leq \left| \int \psi_{B} \; dm_{B,L} -  \int \psi_{B} \; dm_{D,T} \right| + \left| \int \psi_{B} \; dm_{D,T} -  \int \psi_{D} \; dm_{D,T} \right | \\
&\leq  C d((B,L),(D,T))^\alpha  .
\end{align*}\end{proof}

\begin{remark}\label{compute}
The H\"older coefficient $\alpha$ above is computable based on the input data. More precisely, we iterate the cocycle $(A, K)$ a sufficient number $n$ of times, until the inequality~\eqref{compute eq} holds (the existence of such a number of iterates is guaranteed by the our assumptions). Then $\alpha$ is chosen such that $1- \alpha +C \frac{\alpha^2}{2} < 1$, where the constant $C$  depends explicitly on the data. 
\end{remark}

\subsection*{Acknowledgments} A.C. was supported by a FAPERJ postdoctoral grant. M.D. and A.M. were supported by a CNPq doctoral fellowship. S.K. was supported by the CNPq research grant 313777/2020-9 and  by the Coordena\c{c}\~ao de Aperfei\c{c}oamento de Pessoal de N\'ivel Superior - Brasil (CAPES) - Finance Code 001.
\bigskip

\bibliographystyle{amsplain}
\bibliography{references}

\end{document}